\documentclass[12pt]{amsart}
\usepackage{amssymb}
\usepackage{amsthm}
\usepackage{amscd}

\usepackage{amsmath}
\usepackage{epic,eepic}
\usepackage{graphicx}
\DeclareMathAlphabet{\mathpzc}{OT1}{pzc}{m}{it}

\usepackage{color}
\usepackage{latexsym}

\usepackage{pifont}

\theoremstyle{plain}
\newtheorem{lemma}{Lemma}[section]
\newtheorem{prop}[lemma]{Proposition}
\newtheorem{thm}[lemma]{Theorem}
\newtheorem{cor}[lemma]{Corollary}
\newtheorem{aplemma}{Lemma~A.\hspace{-1.5mm}}
\newtheorem{approp}{Proposition~A.\hspace{-1.5mm}}
\newtheorem{apthm}{Theorem~A.\hspace{-1.5mm}}
\newtheorem{apcor}{Corollary~A.\hspace{-1.5mm}}
\newtheorem{intthm}{Theorem}

\usepackage[all]{xy}

\newtheorem{conj}[lemma]{Conjecture}

\theoremstyle{definition}

\newtheorem{rema}[lemma]{Remark}

\newtheorem{remb}{Remark}

\newtheorem{defi}[lemma]{Definition}
\newtheorem{exa}[lemma]{Example}
\newtheorem{aprem}{Remark~A.\hspace{-1.5mm}}
\newtheorem{apdefi}{Definition~A.\hspace{-1.5mm}}
\newcommand{\bde}{\begin{defi}}
\newcommand{\ede}{\end{defi}\vspace{1mm}}
\newcommand{\ble}{\begin{lemma}}
\newcommand{\ele}{\end{lemma}}
\newcommand{\bpr}{\begin{prop}}
\newcommand{\epr}{\end{prop}}
\newcommand{\bt}{\begin{thm}}
\newcommand{\et}{\end{thm}}
\newcommand{\bco}{\begin{cor}}
\newcommand{\eco}{\end{cor}}
\newcommand{\bre}{\begin{rema}}
\newcommand{\ere}{\end{rema}}
\newcommand{\brea}{\begin{rema}}
\newcommand{\erea}{\end{rema}\vspace{1mm}}
\newcommand{\breb}{\begin{remb}}
\newcommand{\ereb}{\end{remb}\vspace{1mm}}
\newcommand{\bex}{\begin{exa}}
\newcommand{\eex}{\end{exa}}
\newcommand{\bpf}{\begin{proof}}
\newcommand{\epf}{\end{proof}\vspace{1mm}}

\newcommand{\bade}{\begin{apdefi}}
\newcommand{\eade}{\end{apdefi}}
\newcommand{\bale}{\begin{aplemma}}
\newcommand{\eale}{\end{aplemma}}
\newcommand{\bapr}{\begin{approp}}
\newcommand{\eapr}{\end{approp}}
\newcommand{\bat}{\begin{apthm}}
\newcommand{\eat}{\end{apthm}}
\newcommand{\baco}{\begin{apcor}}
\newcommand{\eaco}{\end{apcor}}
\newcommand{\bare}{\begin{aprem}}
\newcommand{\eare}{\end{aprem}}

\newcommand{\be}{\begin{enumerate}}
\newcommand{\ee}{\end{enumerate}}
\newcommand{\bcd}{\[\begin{CD}}
\newcommand{\ecd}{\end{CD}\]}
\newcommand{\bit}{\begin{itemize}}
\newcommand{\eit}{\end{itemize}}
\newcommand{\bq}{\begin{quote}}
\newcommand{\eq}{\end{quote}}
\newcommand{\ba}{\begin{array}}
\newcommand{\ea}{\end{array}}

\newcommand{\mcB}{\mathcal{B}}

\newcommand{\mcE}{\mathcal{E}}

\newcommand{\mcL}{\mathcal{L}}

\newcommand{\mcO}{\mathcal{O}}

\newcommand{\mbC}{\mathbb{C}}
\newcommand{\mbD}{\mathbb{D}}

\newcommand{\mbF}{\mathbb{F}}
\newcommand{\mbG}{\mathbb{G}}

\newcommand{\mbP}{\mathbb{P}}

\newcommand{\mbZ}{\mathbb{Z}}

\newcommand{\mfb}{\mathfrak{b}}
\newcommand{\mfc}{\mathfrak{c}}

\newcommand{\mfg}{\mathfrak{g}}
\newcommand{\mfh}{\mathfrak{h}}

\newcommand{\mfl}{\mathfrak{l}}

\newcommand{\mfn}{\mathfrak{n}}

\newcommand{\mfs}{\mathfrak{s}}
\newcommand{\mft}{\mathfrak{t}}

\newcommand{\migi}{\rightarrow}

\newcommand{\isom}{\stackrel{\sim}{\migi}}

\newcommand{\migiincl}{\hookrightarrow}

\newcommand{\migisurj}{\twoheadrightarrow}

\newcommand{\SSP}{\vspace{0mm}}
\newcommand{\LSP}{\vspace{0mm}}

\newcommand{\mr}{\mathrm}
\newcommand{\hidden}[1]{\,}

\begin{document}

\title[Gaudin model modulo $p$, Tango structures, and dormant Miura opers]
{Gaudin model modulo $p$,  Tango structures, \\ and dormant Miura opers}
\author{Yasuhiro Wakabayashi}
\date{}
\markboth{Yasuhiro Wakabayashi}{}
\maketitle
\footnotetext{Y. Wakabayashi: Graduate School of Information Science and Technology, Osaka University, 1-5 Yamadaoka, Suita, Osaka 565-0871, JAPAN;}
\footnotetext{e-mail: {\tt wakabayashi@ist.osaka-u.ac.jp};}
\begin{abstract}
In the present paper, we study the Bethe ansatz equations for Gaudin model and Miura opers in characteristic $p>0$. Our study is based on a work by E. Frenkel, in which solutions to the Bethe ansatz equations are described in terms of Miura opers on the complex projective line. The main result of the present paper provides a positive characteristic analogue of this description. We pay particular attention to the case of Miura $\mathrm{PGL}_2$-opers because dormant generic Miura $\mathrm{PGL}_2$-opers correspond bijectively to Tango structures, which bring various sorts of exotic  phenomena in positive characteristic, e.g., counter-examples to the Kodaira vanishing theorem. As a consequence, we construct new examples of Tango structures by means of solutions to the Bethe ansatz equations modulo $p$.
\end{abstract}
\tableofcontents 

\section*{Introduction}
\SSP

The {\it Bethe ansatz equations}   for a simple finite-dimensional Lie algebra $\mfg$ over  the field of complex numbers $\mbC$ (cf. ~\cite[\S\,3.2, (3.5)]{Fre2}) are the system of  equations on the set of distinct complex numbers $z_1, \cdots, z_m$ ($m \geq 1$) of the form
\begin{align} \label{E003}
\sum_{i=1}^r \frac{\langle \alpha_j, \check{\lambda}_{i}\rangle}{z_j -x_i}  - \sum_{s \neq j} \frac{\langle \alpha_j, \check{\alpha}_{s} \rangle}{z_j -z_s} =0 \hspace{10mm} \left(j =1, \cdots, m\right),
\end{align}
where 
$\alpha_1, \cdots, \alpha_m$ are simple positive roots,
$\check{\alpha}_1, \cdots, \check{\alpha}_m$ are the corresponding coroots, 
$x_1, \cdots, x_r$ ($r\geq 0$) are distinct complex numbers, 
and $\check{\lambda}_1, \cdots, \check{\lambda}_r$ are dominant coweights of $\mfg$.
It is well-known that each solution to  the Bethe ansatz equations (\ref{E003}) 
specifies 
an eigenvector of the Hamiltonian of 
a certain spin model called the {\it Gaudin model}.
It gives an effective way to solve the problem of simultaneous diagonalization of the Gaudin Hamiltonian.

On the other hand, the Bethe Ansatz equations may also be  interpreted as a ``monodromy-free" condition on  Miura opers.
A {\it Miura oper} (cf. \S\,\ref{SS1331}) is, by definition,   an oper equipped with additional data, or more precisely, 
a principal  bundle over an algebraic curve equipped with two Borel reductions and a flat connection satisfying some conditions (including a certain transversality condition).
For example,  $\mr{PGL}_2$-opers and Miura $\mr{PGL}_2$-opers on a complex  algebraic curve 
may be identified with, respectively,  projective and affine structures on the associated Riemann surface.

 In a work by E. Frenkel (cf.  ~\cite{Fre2} and ~\cite{Fre1}),
  a canonical correspondence between  solutions to the Bethe ansatz equations and certain Miura opers with trivial monodromy on the complex projective line
  was  constructed.
This may be regarded as an example of the geometric Langlands correspondence for the projective line.

The present paper discusses  Frenkel's correspondence
{\it in characteristic $p$}, where $p$ is a prime number.
(The previous study of the Bethe ansatz equations in positive characteristic can be found in  ~\cite{Var}.)
Note that we cannot adopt directly the results and their proofs in  \cite{Fre2} or   ~\cite{Fre1}
 because some of them  are based on  properties and techniques  inherent in characteristic $0$.
For example, on a flat connection in characteristic $p$,  the condition   of having trivial monodromy is not sufficient to conclude the local triviality of that  connection. 

Taking this into account, 
 we establish a positive-characteristic version of Frenkel's correspondence
 (cf. Theorem \ref{ThgA} below), asserting that mod $p$ solutions to the Bethe ansatz equations of certain types  
    correspond bijectively to
 generic Miura 
 opers on a pointed projective line  with conditions imposed on their $p$-curvatures.
    (For the previous works concerning such Miura opers, we refer to ~\cite{Wak7} and ~\cite{Wak8}.)

In what follows, we shall describe the main theorems in the present paper.
Let $k$ be an algebraically closed field of characteristic $p$, $G$ a connected simple algebraic group over $k$ of adjoint type satisfying the condition $(*)_G$ introduced in \S\,\ref{SS01}.
Let ${\bf x} := (x_1, \cdots, x_{r+1})$ ($r \geq 0$) be an ordered collection of  $r+1$ distinct closed points of  the projective line $\mbP$ over $k$  with $x_{r+1} = \infty$.
Also,  let $\check{\pmb{\lambda}}:= (\check{\lambda}_1, \cdots, \check{\lambda}_{r+1})$ be an element of $(\mft_\mr{dom}^F)^{r+1}$ (cf.  (\ref{Ee121}) for the definition of $\mft_\mr{dom}^F$) and  $\pmb{\alpha} := (\alpha_1, \cdots, \alpha_m)$ an element of $\Gamma^m$, where  $m \geq 1$ and $\Gamma$ denotes the set of simple positive roots with respect to a fixed pair of a  maximal torus and a Borel subgroup of $G$.
 Suppose that these collections of data  satisfy the following equality
 \begin{align} \label{Ej6789}
 -\sum_{i=1}^{r+1} \check{\lambda}_i + \sum_{j=1}^m \check{\alpha}_j = 2 \check{\rho}
 \end{align}
(cf.  (\ref{W100})), where $\check{\rho}$ denotes the sum of the fundamental coweights of simple positive roots. 

 Denote by $C^m ({\bf x})$  (cf. (\ref{Ee449})) the set of ordered collections of $m$ distinct closed points in $\mbP \setminus \{x_1, \cdots, x_{r+1}\}$.
Also, denote by
\begin{align}
G\text{-}\mr{BA}_{\check{\pmb{\lambda}}, \pmb{\alpha}} 
\end{align}
(cf. (\ref{Ef22})) the subset of $C^m ({\bf x})$ consisting of elements ${\bf z}:= (z_1, \cdots, z_m)$ satisfying the Bethe ansatz equations (\ref{E003}) (considered as a system of equations with coefficients in $k$) associated to  the fixed triple $({\bf x}, \check{\pmb{\lambda}}, \pmb{\alpha})$.
On the other hand, we have a certain  set
\begin{align}
G \text{-} \mr{MOp} (\mbP^{D_{{\bf x}}\text{-}\mr{log}}; -\check{\pmb{\lambda}} -\check{\rho})_{\mr{triv}, + \pmb{\alpha}^W}^{^{p\text{-}\mr{nilp}}}
\  \left(\text{resp.,} \ G \text{-} \mr{MOp} (\mbP^{D_{{\bf x}}\text{-}\mr{log}}; -\check{\pmb{\lambda}} -\check{\rho})_{\mr{gen}, + \pmb{\alpha}^W}^{^\mr{Zzz...}} \right)
\end{align}
consisting of  $p$-nilpotent generic  Miura $G$-opers  satisfying a monodromy-free condition (resp., dormant generic $G$-opers)  on $\mbP$ equipped  with a log structure;
 the precise definition of this set can be seen in (\ref{Ee556}) (resp.,  (\ref{Ee600})).  
Since any connection with  vanishing $p$-curvature is locally trivial,  the  dormancy condition is stronger than the monodromy-free condition.
This implies that there exists a canonical inclusion 
 \begin{align} \label{Ee409}
G \text{-} \mr{MOp} (\mbP^{D_{{\bf x}}\text{-}\mr{log}}; -\check{\pmb{\lambda}} -\check{\rho})_{\mr{gen}, + \pmb{\alpha}^W}^{^\mr{Zzz...}} \migiincl  G \text{-} \mr{MOp} (\mbP^{D_{{\bf x}}\text{-}\mr{log}}; -\check{\pmb{\lambda}} -\check{\rho})_{\mr{triv}, + \pmb{\alpha}^W}^{^{p\text{-}\mr{nilp}}}
 \end{align}
 (cf.  Propositions \ref{P015}, (i), and \ref{pgT011}).

Then, the  following theorem
 tells us relationships between these sets of Miura $G$-opers and the solutions to the Bethe ansatz equations modulo $p$.

\vspace{3mm}
\begin{intthm}[cf. Theorem \ref{TT011}, Proposition \ref{pgT011}] \label{ThgA} 
\begin{itemize}
\item[(i)]
There exists a canonical  bijection  of sets
 \begin{align} \label{Ee55}
 G \text{-} \mr{BA}_{\check{\pmb{\lambda}}, \pmb{\alpha}} 
 \isom
 G \text{-} \mr{MOp} (\mbP^{D_{{\bf x}}\text{-}\mr{log}}; -\check{\pmb{\lambda}} -\check{\rho})_{\mr{triv}, + \pmb{\alpha}^W}^{^{p\text{-}\mr{nilp}}}.
 \end{align}
  \item[(ii)]
 If, $\check{\pmb{\lambda}}= 0^{r+1}  \left(:= (0, 0, \cdots, 0) \right)$, then the injection (\ref{Ee409}) becomes bijective. In particular, 
  we obtain a bijection of sets
 \begin{align}
  G \text{-} \mr{BA}_{0^{r+1}, \pmb{\alpha}} 
 \isom
 G \text{-} \mr{MOp} (\mbP^{D_{{\bf x}}\text{-}\mr{log}};  (-\check{\rho})^{r+1})_{\mr{gen}, + \pmb{\alpha}^W}^{^\mr{Zzz...}}.
 \end{align}
 \end{itemize}
 \end{intthm}


The above theorem for $G = \mr{PGL}_2$ provides 
an application to  the problem of constructing (pre-)Tango structures. 
A {\it (pre-)Tango structure}  is  a  certain line bundle on an algebraic curve (cf. Definition  \ref{W400}, (i),  or Definition \ref{Ww551} for the precise definition) and  has an important feature in that
it brings various sorts of exotic phenomena in positive characteristic.
For example,  each Tango structure on a proper smooth curve  in characteristic $p$ associates an algebraic surface violating the Kodaira vanishing theorem, as well as other reasonable theorems valid  in characteristic $0$ (cf. ~\cite[\S\,7.2]{Wak7} for the construction of such surfaces using  Tango structures).  
To the author's knowledge, 
there are   very few previous  examples of Tango structures  described explicitly.
In fact, the existence of a Tango structure implies a strong restriction on the genus $g$  of the underlying curve, i.e., $p$ must divide $2g-2$.

However, by combining  Theorem \ref{ThgA} with  the correspondence between Tango structures and dormant generic Miura $\mr{PGL}_2$-opers  proved in ~\cite{Wak7}, 
we  can construct (cf. Theorem \ref{ThB} below) infinitely many examples of Tango curves, i.e., algebraic  curves admitting a Tango structure (cf. Definition \ref{W400},  (ii)).
The well-known examples  given 
by M. Raynaud  (cf. ~\cite[Example]{Ray} or ~\cite[Example 1.3]{Muk})
 may be thought of as  special cases
  of our  construction, and  we obtain   
  other new  Tango curves.
 

\SSP
\begin{intthm}[cf. Theorem \ref{T041}] \label{ThB} 
Suppose that $r=0$ and that we are given 
a pair of positive integers   $(a, b)$
 with $\mr{gcd}(a, bp-1) =1$ and $m = ap$.
Let   $(z_1, \cdots, z_{ap})$ be a solution to the Bethe ansatz equation classified by $\mr{PGL}_2 \text{-}\mr{BA}_{0, \pmb{\alpha}}$.
Denote by $Y$ 
   the desingularization of the plane curve defined by the equation
\begin{align}
y^{bp -1} = \prod_{j=1}^{ap} (x -z_j).
\end{align}
Then, $Y$
is a Tango curve.
 \end{intthm}

\vspace{10mm}
\section{Opers and Miura opers} 
\SSP

In this section, we recall the definitions of an oper and a (generic) Miura oper.

\LSP
\subsection{Algebraic groups and Lie algebras} \label{SS01}

Throughout the present paper, let us fix a prime $p$ and  an algebraically closed field $k$ of characteristic $p$.
Also, let $G$ be 
  a connected simple algebraic group  over $k$ of adjoint type satisfying 
the 
condition $(*)_G$ described as follows.
\begin{itemize}
\item[$(*)_G$ :]
$G$ is either equal to $\mr{PGL}_n$ with $1 <n <p$ or satisfies the inequality  $p> 2 \cdot h$, where $h$ denotes the Coxeter number of $G$.
\end{itemize}

Let us fix  a maximal torus $T$
  of $G$ and 
  a Borel subgroup $B$
    of $G$ containing $T$.
Write $N = [B, B]$, i.e., the unipotent radical of $B$, and write $W$  for  the Weyl group of $(G, T)$.
Denote by $\mfg$, $\mfb$, $\mfn$, and  $\mft$   the  Lie algebras of $G$, $B$,  $N$,  and $T$, respectively (hence $\mft,  \mfn \subseteq \mfb \subseteq \mfg$).
Denote by 
$\Phi^+$  the set of positive  roots in $B$ with respect to $T$ and by   $\Phi^-$ the set of negative roots. 
Also, denote by $\Gamma$ ($\subseteq \Phi^+$) the set of  simple positive  roots.
Given each character $\varphi : T \migi \mbG_m$ (resp., each cocharacter $\check{\varphi} : \mbG_m \migi T$), we shall  use, by abuse of notation, the same notation $\varphi$ (resp., $\check{\varphi}$) to denote  its differential $d \varphi \in \mft^\vee$ (resp., $d\check{\varphi} \in \mft$). 

For each $\alpha \in \Phi^+ \cup \Phi^-$, we write
\begin{equation} \label{22223}
\mfg^\alpha := \big\{ x \in \mfg \ \big| \ \text{$\mr{ad}(t)(x) =  \alpha  (t) \cdot x$ for all  $t\in T$}  \big\}.\end{equation}
Each $\mfg^{-\alpha}$ ($\alpha \in \Gamma$)
may be thought of as a subspace of $\mfg / \mfb$ closed under the adjoint $B$-action.
The Lie algebra $\mfg$ is equipped with the principal gradation $\mfg = \bigoplus_{i = - \mr{rk}(\mfg)}^{\mr{rk}(\mfg)} \mfg_i$,  which
restricts to  identifications $\mft = \mfg_0$, $\bigoplus_{\alpha \in \Gamma} \mfg^\alpha = \mfg_1$,  $\bigoplus_{\alpha \in \Gamma} \mfg^{-\alpha} = \mfg_{-1}$, and 
 $\mfn = \bigoplus_{i =1}^{\mr{rk}(\mfg)}\mfg_i$.

We shall write
$\check{\rho}$ for the 
cocharacter
 $\mbG_m \migi T$ (as well as its differential)
   defined as   the sum $\sum_{\alpha \in \Gamma}  \check{\omega}_\alpha$ such that  each $\check{\omega}_\alpha$ ($\alpha \in \Gamma$) denotes  the fundamental coweight  of $\alpha$.
We  fix a generator $f_\alpha$ of $\mfg^{-\alpha}$  for each $\alpha \in \Gamma$.
Hence, we obtain $p_{-1} := \sum_{\alpha \in \Gamma} f_\alpha$.

If $\mft_\mr{reg}$ denotes the set of regular elements in $\mft$, then
it follows from ~\cite[Chap.\,VI,  Theorem 7.2]{KW}  that 
\begin{align} \label{Ee120}
\mft_\mr{reg} = \left\{ \check{\lambda} \in \mft   \ | \ \alpha (\check{\lambda}) \neq 0 \ \text{for any root $\alpha \in \Phi^+ \cup \Phi^{-}$}  \right\}.
\end{align}
We obtain   the set of $F$-invariant elements $\mft^F_\mr{reg}$  in $\mft_\mr{reg}$, where $F$ denotes the  Frobenius endomorphism of $\mft_\mr{reg}$ viewed as a $k$-scheme, and
write
\begin{align} \label{Ee121}
\mft_\mr{dom}^F :=  \left\{ \check{\lambda} \in \mft \, \big| \,  \check{\lambda} + \check{\rho} \in \mft_\mr{red}^F  \right\}.
\end{align}

\LSP
\subsection{Opers} \label{SS02}

Let $X$ be a connected proper smooth
  curve over $k$ and $D$
 a reduced effective divisor on $X$.
One can equip  $X$  with  a log structure  determined by $D$ in the usual manner;
we denote the resulting log scheme by $X^{D \text{-} \mr{log}}$.
In particular, if $D=0$, then $X^{D \text{-} \mr{log}}$ coincides with $X$.

Let
$G_0$ be an algebraic group over $k$ and 
 $\mcE$    a (right) $G_0$-bundle on  $X$.
Given a $k$-vector space $\mfh$ equipped with a (left) $G_0$-action,
 we shall write $\mfh_{\mcE}$ for the vector bundle on $X$ associated with the relative affine space $\mcE \times^{G_0} \mfh$ 
($:= (\mcE \times_k \mfh) /G_0$) over $X$.
By  a {\bf $D$-log connection} on $\mcE$, we mean a logarithmic   connection on $\mcE$ with respect to the log structure of $X^{D \text{-}\mr{log}}$.

Let $G$ be  as in the previous subsection and 
$\mcE$  a  $G$-bundle on  $X$.
The adjoint $G$-action  on $\mfg$ gives rise to
 a vector bundle $\mfg_{\mcE}$, i.e., the adjoint vector bundle associated to $\mcE$.

Next, suppose that we are given
a $D$-log 
 connection $\nabla$  on   $\mcE$ and  a $B$-reduction $\mcE_B$  of $\mcE$ (i.e., a $B$-bundle $\mcE_B$ on $X$ together with an isomorphism of $G$-bundles  $\mcE_B \times^B G \isom \mcE$).
 Let us  choose, locally on $X$,  a $D$-log connection $\nabla'$ on $\mcE$ preserving $\mcE_B$, and take the difference $\nabla - \nabla'$, which specifies   a section of $\Omega_{X^{D \text{-} \mr{log}}/k} \otimes \mfg_{\mcE_B} \left(=\Omega_{X^{D \text{-} \mr{log}}/k}\otimes \mfg_{\mcE}\right)$.
The local section of  $\Omega_{X^{D \text{-}\mr{log}}/k} \otimes (\mfg/\mfb)_{\mcE_B}$ determined by this section via projection
does not depend on the choice of  $\nabla'$.
Hence, these sections defined for various $\nabla'$'s 
may be glued together to obtain  a global section of $\Omega_{X^{D \text{-} \mr{log}}/k} \otimes (\mfg/\mfb)_{\mcE_B}$; we shall denote this section
 by $\nabla /\mcE_B$.

A {\bf $G$-oper} on 
$X^{D \text{-} \mr{log}}$ means
 a triple 
$\mcE^\spadesuit := (\mcE, \nabla, \mcE_B)$, where
$\mcE$, $\nabla$, and $\mcE_B$ are as above such that 
the section  $\nabla/ \mcE_B$ lies in the submodule $\Omega_{X^{D \text{-} \mr{log}}/k} \otimes (\bigoplus_{\alpha \in \Gamma} \mfg^{-\alpha})_{\mcE_B}$
 and that  its  image in $\Omega_{X^{D \text{-} \mr{log}}/k} \otimes \mfg^{- \beta}_{\mcE_B}$ (for each $\beta \in \Gamma$) via the projection $\bigoplus_{\alpha \in \Gamma} \mfg^{-\alpha} \migisurj \mfg^{-\beta}$ specifies 
a nowhere vanishing section.
In a natural manner, one can define the notion of an isomorphism between $G$-opers.
Thus, we obtain the set
\begin{align}
G \text{-} \mr{Op}(X^{D \text{-} \mr{log}})
\end{align}
of isomorphism classes of $G$-opers on $X^{D \text{-}\mr{log}}$.

Next,
let $\widehat{\mcO}$ be a complete discrete valuation ring  over $k$ whose residue field is isomorphic to $k$, and 
write $\mbD$ for the formal disc $\mr{Spec}(\widehat{\mcO})$.
The closed point of $\mbD$, denoted by $x_0$,   defines a reduced effective divisor, and hence, defines a log structure on $\mbD$; we denote the resulting log scheme by $\mbD^\mr{log}$.
Also,  we  write $\mbD^\times := \mbD \setminus \{ x_0 \}$.
As in the case of the entire curve 
  $X^{D \text{-}\mr{log}}$ discussed above, we have  the definition of a $G$-oper on $\mbD$
(resp., $\mbD^\mr{log}$; resp., $\mbD^\times$).
Denote by
\begin{align}
G \text{-}\mr{Op} (\mbD)
 \left(\text{resp.,} \ G \text{-}\mr{Op} (\mbD^\times)\right)
\end{align}
the set of isomorphism classes of $G$-opers on $\mbD$ (resp.,    $\mbD^\times$).

 Let $\mcE^\spadesuit := (\mcE, \nabla, \mcE_B)$ be a $G$-oper on $\mbD^\mr{log}$.
 After choosing a uniformizer $t$ of $\widehat{\mcO}$ (which gives an isomorphism $k[\![t]\!] \isom \widehat{\mcO}$) and applying a suitable gauge transformation,
 $\nabla$ may be expressed as a log connection on the trivial $G$-bundle $\mbD \times G$ of the form
 \begin{align} \label{Ee800}
 \partial_t + \frac{1}{t} (p_{-1} + u) + {\bf u}(t)
 \end{align}
 for some $u \in \mfb$, ${\bf u}(t) \in \mfb [\![t]\!] := \mfb (k[\![t]\!])$, where $\partial_t := \frac{d}{dt}$ (cf. ~\cite[Definition 1.20, Proposition 2.8]{Wak5}).
In particular, the monodromy operator,   in the sense of ~\cite[Definition 1.46]{Wak5},  of $\nabla$ at $x_0$   (with respect to this expression)  is $p_{-1}+u$.
 Let 
 $\mfc  \left(:= \mfg \ooalign{$/$ \cr $\,/$}\hspace{-0.5mm}G\right)$
  denote the GIT quotient of  $\mfg$ by 
the adjoint $G$-action,  and let 
 $\chi : \mfg \migi \mfc$
  denote the natural quotient.
The element  $\rho := \chi (p_{-1}+ u) \in \mfc (k)$ depends neither  on  the choice of $t$ nor the expression (\ref{Ee800}) of $\nabla$;
  we  say that $\mcE^\spadesuit$ is {\bf of radius $\rho$}  (cf. ~\cite[Definition 2.29]{Wak5}).

 For each $\rho \in \mfc (k)$, we shall write
 \begin{align}
 G \text{-}\mr{Op} (\mbD^\mr{log}; \rho)
 \end{align}
 for the
 set of isomorphism classes of $G$-opers on $\mbD^\mr{log}$
   of radius $\rho$.
Note that the following assertion was already proved in    ~\cite[Proposition 2.1.1  and the discussion at the end of  \S\,9.1]{Fr} when the base field is replaced with $\mbC$.
 But, since the case of positive characteristic requires special consideration on the existence of the exponential map $\mfn \migi N$, we here reprove it.

\SSP
\bpr\label{P395f} The  maps of sets
 \begin{align} \label{Ej2}
 G \text{-}\mr{Op}(\mbD) \migi G \text{-}\mr{Op}(\mbD^\times), \hspace{10mm} 
 G \text{-} \mr{Op} (\mbD^\mr{log}; \chi (- \check{\rho})) \migi 
 G \text{-}\mr{Op}(\mbD^\times)
 \end{align}
 defined by restriction via the inclusions $\mbD^\times \migiincl \mbD$ and $\mbD^\times \migiincl \mbD^\mr{log}$, respectively, are injective.
 Moreover, the image of the former map is contained in that of the latter map,
  and the resulting injection 
  \begin{align} \label{Ee8801}
 G \text{-}\mr{Op}(\mbD) \migiincl G \text{-} \mr{Op} (\mbD^\mr{log}; \chi (- \check{\rho}))
 \end{align}
can be given by gauge transformation by 
$\check{\rho} \circ t \in T(k(\!(t)\!))$ (for any uniformizer $t \in \widehat{\mcO}$).
  \epr
\begin{proof} 
To prove the former assertion, we only consider the injectivity of 
$ G \text{-}\mr{Op}(\mbD) \migi G \text{-}\mr{Op}(\mbD^\times)$ because
the proof of the remaining one is entirely similar.
Let $\mcE^\spadesuit_1$ and $\mcE^\spadesuit_2$ be $G$-opers on $\mbD$  such that  $\mcE^\spadesuit_1 \cong \mcE^\spadesuit_2$ when restricted to $\mbD^\times$.
After choosing a uniformizar $t$ of $\widehat{\mcO}$ and apply a suitable gauge transformation, we can express  the connection $\nabla_l$  defining  $\mcE^\spadesuit_i$ (for each $l=1,2$)   as 
\begin{align}
\nabla_l = \partial_t + p_{-1} + {\bf u}_l (t)
\end{align}
for some   ${\bf u}_l (t) \in \mfb [\![t]\!]$.
By assumption,  there exists an element $b \in B (k(\!(t)\!))$ satisfying the equality 
$\nabla_2 = (\nabla_1)_b$, where $(\nabla_1)_b$ denotes the connection obtained from $\nabla_1$ by carrying out    the gauge transformation by $b$.
The problem is reduced to proving   that $b \in B(k[\![t]\!])$.
Suppose, on the contrary, that $b \notin B(k[\![t]\!])$.
Since the  entry lying in $\mfg_{-1}  \left(\subseteq \mfg \right)$ of $\nabla_1$ 
is the same as that of $\nabla_2$
 (i.e., coincides with $p_{-1}$),
$b$ turns out to be an element of $N (k(\!(t)\!))$. 
By the assumption $(*)_G$, we have the exponential map $\mr{exp} : \mfn \isom N$ given by ~\cite[Proposition 1.31 and Remark 1.35]{Wak5}.
In particular, $b = \mr{exp}({\bf v}(t))^{-1}$ for some ${\bf v}(t) \in \mfn (\!(t)\!)$.
It follows from   ~\cite[Corollary 1.34]{Wak5} that the equality $\nabla_2 = (\nabla_1)_b$ implies 
\begin{align} \label{Ej10}
\partial_t + p_{-1} + {\bf u}_2 (t) = \partial_t + \frac{d}{dt} (-{\bf v}(t)) + \sum_{s=0}^\infty \frac{1}{s!} \cdot \mr{ad}({\bf v}(t))^s (p_{-1}+ {\bf u}_1(t)). 
\end{align}
According to the principal  gradation on $\mfg$,  the elements ${\bf u}_1(t)$, ${\bf u}_2 (t)$, and ${\bf v}(t)$ can be decomposed as 
\begin{align}
{\bf u}_1(t) = \sum_{i=1}^{\mr{rk}(\mfg)} {\bf u}_{1, i} (t),
 \ \ {\bf u}_2(t) = \sum_{i=1}^{\mr{rk}(\mfg)} {\bf u}_{2, i} (t),
 \ \  \text{and} \ \   
 {\bf v}(t) = \sum_{i=1}^{\mr{rk}(\mfg)} {\bf v}_i (t),
 \end{align}
  respectively, where ${\bf u}_{1, i} (t), {\bf u}_{2, i} (t) \in \mfg_i [\![t]\!]$ and ${\bf v}_i (t)\in \mfg_i (\!(t)\!)$ ($i = 1, \cdots, \mr{rk}(\mfg)$).
Since $b \notin B (k[\![t]\!])$ (or equivalently, ${\bf v}(t) \notin \mfn (k[\![t]\!])$), 
  the positive integer 
  \begin{align}
  i_\mr{min} :=\mr{min} \left\{ i \in \mbZ_{>0} \, | \, {\bf  v}_i (t) \notin \mfg_i [\![t]\!]\right\}
  \end{align}
   is well-defined.
Given a pair of positive integers $(i, j)$ with  $j \leq i \leq   \mr{rk}(\mfg)$, 
we denote by $F_j^i$ the degree $i$ graded $k$-linear endomorphism  of $\mfg (\!(t)\!)$ defined as
\begin{align} \label{Ej12}
F_j^i := \sum_{s=j}^i \frac{1}{s!} \cdot \sum_{\substack{(l_1, \cdots, l_s) \in \mbZ_{>0}^s, \\ l_1 + \cdots + l_s = i}} \mr{ad} ({\bf v}_{l_1}(t)) \circ \cdots \circ \mr{ad}({\bf v}_{l_s}(t)) : \mfg (\!(t)\!) \migi \mfg (\!(t)\!).
\end{align}
By comparing the respective entries lying in $\mfg_{i_\mr{min}-1}(k(\!(t)\!))$ of the both sides of (\ref{Ej10}), we obtain the equality
\begin{align}
{\bf u}_{2, i_\mr{min}-1} =  [{\bf v}_{i_\mr{min}}(t), p_{-1}] + F_2^{i_\mr{min}}(p_{-1}) + \sum_{i=0}^{i_\mr{min}-1} F_{1}^{i_\mr{min}-1-i} ({\bf u}_{1, i}(t)).
\end{align}
It follows from the definition of $i_\mr{min}$ that
this equality implies $[{\bf v}_{i_\mr{min}}(t), p_{-1}] \in \mfg_{i_\mr{max} -1}[\![t]\!]$.
Moreover, since 
 $[-, p_{-1}]  \left(= - \mr{ad} (p_{-1}) \right) : \mfg_{i_\mr{min}} \migi \mfg_{i_\mr{min}-1}$ is injective,  we  have ${\bf v}_{i_\mr{min}}(t) \in \mfg_{i_\mr{min}}[\![t]\!]$.
 This is a contradiction, and hence, completes the proof of the former assertion.

The latter assertion follows from the definitions of various maps involved, so we finish the proof of the proposition.
\end{proof} 
\SSP

\begin{rema}\label{Rr90}
 Since the condition  $(*)_G$ is assumed to be fulfilled,
  it follows from Chevalley's theorem (cf. ~\cite[Theorem 1.1.1]{Ngo},  ~\cite[Chap.\,VI, Theorem 8.2]{KW}) that
 the composite $\mft \migiincl \mfg \stackrel{\chi}{\migi} \mfc$ induces an isomorphism 
  $\mft/W \isom  \mfc$.
 In particular, an element $\check{\lambda} \in \mft$ satisfies the equality $\chi (\check{\lambda}) = \chi (-\check{\rho})$ if and only if $\check{\lambda} = w (-\check{\rho})$ for some $w \in W$.
 \end{rema}
\SSP

Let $X$  be as before and ${\bf x} := (x_i)_{i=1}^r$ ($r \geq 1$)  an ordered collection of distinct closed points of $X$.
 We shall write $D_{{\bf x}}:= \sum_{i=1}^r [x_i]$, where $[x_i]$ ($i=1, \cdots, r$) denotes the reduced effective  divisor on $X$ determined by $x_i$.
 For each $i \in \{1, \cdots, r \}$, denote by $\mbD_{x_i}$ the formal neighborhood of $x_i$ in $X$.
 Given an element   $\pmb{\rho} := (\rho_i)_{i=1}^r$ of $\mfc (k)^{r}$, 
 we say that a $G$-oper on $X^{D_{{\bf x}}\text{-} \mr{log}}$ is {\bf of radii $\pmb{\rho}$} (cf. ~\cite[Definition 2.32]{Wak5})  if,  for each $i \in \{1, \cdots, r\}$,  the $G$-oper on  $\mbD_{x_i}^\mr{log}$  defined as    its  restriction      is of radius $\rho_i$.
 Denote by 
 \begin{align}
 G \text{-} \mr{Op} (X^{D_{{\bf x}}\text{-} \mr{log}}; \pmb{\rho})
 \end{align}
 the subset of $G \text{-} \mr{Op} (X^{D_{{\bf x}}\text{-} \mr{log}})$ classifying 
 $G$-opers of radii $\pmb{\rho}$.

\LSP
\subsection{Generic Miura opers} \label{SS1331}

Next, we  recall the notion of  a (generic) Miura oper on a  curve with log structure (cf. ~\cite[Definition 3.2.1]{Wak7}).

A {\bf Miura $G$-oper} on $X$ (resp., on $X^{D \text{-}\mr{log}}$) is a quadruple $\widehat{\mcE}^\spadesuit := (\mcE, \nabla, \mcE_B, \mcE'_B)$, where $(\mcE,  \nabla, \mcE_B)$ is a $G$-oper  on $X$ (resp., on $X^{D \text{-}\mr{log}}$) and $\mcE'_B$ is another  $B$-reduction of $\mcE$   horizontal with respect to $\nabla$.
The definition of an isomorphism between Miura $G$-opers can be formulated in a natural fashion.

Let $\widehat{\mcE}^\spadesuit := (\mcE, \nabla, \mcE_B, \mcE'_B)$ be a Miura $G$-oper on 
either $X$ or  $X^{D \text{-} \mr{log}}$.
By twisting the flag variety $G/B$ by the $B$-bundle  $\mcE_B$,
 we obtain a proper scheme 
$(G/B)_{\mcE_B} := \mcE_B \times^B (G/B)$ over $X$.
Here,  denote by $w_0$ the longest element of $W$.
Note that the Bruhat decomposition $G = \coprod_{w \in W} B  w  B$ gives rise to
 a decomposition
\begin{align}
(G/B)_{\mcE_B} = \coprod_{w \in W} S_{\mcE_B, w}
\end{align}
of $(G/B)_{\mcE_B}$, where
each $S_{\mcE_B, w}$ denotes  the $\mcE_B$-twist of $B w_0 w B$, i.e., $S_{\mcE_B, w} := \mcE_B \times^B (B w_0 w B)$.
The $B$-reduction $\mcE'_B$ determines a section 
$\sigma_{\mcE_B, \mcE'_B} :  X \migi (G/B)_{\mcE_B}$
  of the natural projection $(G/B)_{\mcE_B} \migi X$.
Given an element $w$ of $W$ and a point $x$ of $X$, 
we  say that $\mcE_B$ and $\mcE'_B$ are in {\bf relative position $w$ at $x$} if $\sigma_{\mcE_B, \mcE'_B} (x)$ belongs to $S_{\mcE_B, w}$.
In particular, if $\sigma_{\mcE_B, \mcE'_B} (x)$ belongs to $S_{\mcE_B, 1}$, then we say that $\mcE_{B}$ and $\mcE'_B$ are in {\bf generic position at $x$}.
Moreover,  a  Miura $G$-oper $\widehat{\mcE}^\spadesuit := (\mcE, \nabla, \mcE_B, \mcE'_B)$ 
 is called {\bf generic} if
$\mcE_B$ and $\mcE'_B$ are in generic position at every  point of $X$. 

We also have  the various notions just recalled  in the case where the underlying curve (i.e.,  $X$ or $X^{D \text{-} \mr{log}}$) is replaced with either $\mbD$, $\mbD^\times$, or  $\mbD^{\mr{log}}$.

\LSP
\subsection{Exponents of  Miura opers} \label{SS171}

Let $\widehat{\mcE}^\spadesuit := (\mcE, \nabla, \mcE_B, \mcE'_B)$ be a generic Miura $G$-oper on $\mbD^\mr{log}$.
After  choosing a uniformizer $t$ of $\widehat{\mcO}$ and  applying a gauge transformation, $\nabla$ may be expressed as a log connection on the trivial $G$-bundle $\mbD \times G$  of the form
\begin{align} \label{Ee445}
\partial_t + \frac{1}{t} \left(p_{-1} + \check{\lambda}\right) + {\bf u}(t)
\end{align}
for some
$\check{\lambda}\in\mft$, 
${\bf u} (t)\in\mft [\![t]\!]$
 (cf. ~\cite[Proposition 3.4.3]{Wak7} for the case where the underlying curve is   globally defined).
The element   $\check{\lambda}$ depends neither on the choice of $t$ nor on  the expression (\ref{Ee445}) of $\nabla$.
  In this situation, we say that $\widehat{\mcE}^\spadesuit$ is {\bf of exponent $\check{\lambda}$}.

\SSP
\begin{rema}\label{Rr1019}
Let $\check{\lambda}$ and $w$ be  elements of $\mft_\mr{reg}$ and  $W$, respectively.
Then,  it is verified  that   $\chi (w(\check{\lambda})) = \chi (\check{\lambda}) = \chi (p_{-1}+\check{\lambda})$  (cf. Remark \ref{Rr90}).
In particular, 
if 
$\widehat{\mcE}^\spadesuit := (\mcE, \nabla, \mcE_B, \mcE'_B)$ is  a generic Miura $G$-oper on $\mbD^\mr{log}$ of exponent $\check{\lambda} \left(\in \mft\right)$, then the radius of 
 its  underlying  $G$-oper $(\mcE, \nabla, \mcE_B)$  coincides with   $\chi (\check{\lambda}) \left(= \chi (w(\check{\lambda})) = \chi (p_{-1}+\check{\lambda})\right)$.
  \end{rema}
\SSP

\begin{rema}\label{R1019}
 By applying the gauge transformation by $\check{\rho} \circ t \in T(k(\!(t)\!))$,
we obtain  a bijection from  the set of isomorphism classes of Miura $G$-opers on $\mbD$ to  the set of isomorphism classes of Miura $G$-opers on $\mbD^\mr{log}$ of exponent $-\check{\rho}$.
This bijection is compatible (under the equality $\chi (-\check{\rho}) = \chi (p_{-1}- \check{\rho})$) with the injection (\ref{Ee8801}) via  forgetting the data of the second $B$-reductions defining  Miura $G$-opers.
  \end{rema}
\SSP

Next, let $X$ be as before and  let ${\bf x} := (x_i)_{i=1}^r$,
 $D_{{\bf x}}$, and $\mbD_{x_i}$'s  be as in \S\,\ref{SS02}.
Also,   let $\check{\pmb{\lambda}} := (\check{\lambda}_i)_{i=1}^r \in (\mft_\mr{reg})^{r}$.
  We shall say that a  Miura $G$-oper  $\widehat{\mcE}^\spadesuit := (\mcE, \nabla, \mcE_B, \mcE'_B)$ on $X^{D_{\bf{x}} \text{-} \mr{log}}$  is {\bf of exponents $\check{\pmb{\lambda}}$} 
  if,
  for each $i \in \{1, \cdots, r\}$, 
 the Miura $G$-oper on $\mbD_{x_i}^\mr{log}$ 
 induced by restricting  $\widehat{\mcE}^\spadesuit$
    is generic and of exponent $\check{\lambda}_i$.
  Denote by
  \begin{align} \label{W53}
  G \text{-}\mr{MOp}(X^{D_{{\bf x}}\text{-} \mr{log}}; \check{\pmb{\lambda}})_\mr{gen}
  \end{align}
  the set of isomorphism classes of generic Miura $G$-opers  on $X^{D_{{\bf x}} \text{-} \mr{log}}$ of 
  exponents $\check{\pmb{\lambda}}$.
 
 For each positive integer  $m$, we shall write 
  \begin{align} \label{Ee449}
  C^m ({\bf x})
  \end{align}
   for the set of ordered collections of $m$ distinct closed points in $X \setminus \{x_1, \cdots, x_r \}$. 
Given  each ${\bf w} := (w_1, \cdots, w_m) \in W^{m}$,  we shall set 
  \begin{align} \label{W4}
    G \text{-}\mr{MOp}(X^{D_{{\bf x}}\text{-} \mr{log}}; \check{\pmb{\lambda}})_{\mr{gen}, + {\bf w}} &:= \coprod_{{\bf z} \in C^m ({\bf x})}
       G \text{-}\mr{MOp}(X^{D_{({\bf x}, {\bf z})} \text{-} \mr{log}}; (\check{\pmb{\lambda}}, {\bf w}(-\check{\rho})))_\mr{gen}, 
    \end{align}
    where ${\bf w}(-\check{\rho}):= (w_j (-\check{\rho}))_{j=1}^m \in \mft^{m}$.

 For  each ${\bf z} \in C^m ({\bf x})$, the  injection   (\ref{Ee8801}) resulting from Proposition \ref{P395f} induces 
 an injection
\begin{align} \label{Ee900}
 G \text{-}\mr{Op}(X^{D_{{\bf x}}\text{-} \mr{log}}; \chi (\check{\pmb{\lambda}}))
 \migiincl 
  G \text{-}\mr{Op}(X^{D_{({\bf x}, {\bf z})}\text{-} \mr{log}}; (\chi (\check{\pmb{\lambda}}), \chi (-\check{\rho})^{m})),
\end{align}
 where $ \chi (\check{\pmb{\lambda}}) := (\chi (\check{\lambda}_1), \cdots, \chi (\check{\lambda}_r)) \in \mfc (k)^{r}$ and $\chi (-\check{\rho})^{m}:= (\chi (-\check{\rho}), \cdots, \chi (-\check{\rho}))  \in \mfc (k)^{m}$.
 For  each ${\bf w}$ as above,  we denote by
  \begin{align}\label{Ee805}
   G \text{-}\mr{MOp}(X^{D_{{\bf x}}\text{-} \mr{log}}; \check{\pmb{\lambda}})_{\mr{triv}, + ({\bf z}; {\bf w})}
    \left( \subseteq G \text{-}\mr{MOp}(X^{D_{({\bf x}, {\bf z})}\text{-} \mr{log}}; (\check{\pmb{\lambda}}, {\bf w}(-\check{\rho})))_\mr{gen} \right)
  \end{align}
the inverse image of the image of (\ref{Ee900}) via the forgetting map
\begin{align} \label{W1}
G \text{-}\mr{MOp}(X^{D_{({\bf x}, {\bf z})}\text{-} \mr{log}}; (\check{\pmb{\lambda}}, {\bf w}(-\check{\rho})))_\mr{gen}
\migi
G \text{-}\mr{Op}(X^{D_{({\bf x}, {\bf z})}\text{-} \mr{log}}; (\chi (\check{\pmb{\lambda}}), \chi (-\check{\rho})^{m}))
\end{align}
(cf. Remark
 \ref{Rr1019}).
That is to say, the following square diagram is commutative and cartesian:
\begin{align} \label{W2}
\begin{CD}
   G \text{-}\mr{MOp}(X^{D_{{\bf x}}\text{-} \mr{log}}; \check{\pmb{\lambda}})_{\mr{triv}, + ({\bf z}; {\bf w})} @> \mr{incl.}>> G \text{-}\mr{MOp}(X^{D_{({\bf x}, {\bf z})}\text{-} \mr{log}}; (\check{\pmb{\lambda}}, {\bf w}(-\check{\rho})))_\mr{gen}
\\
@VVV @VV (\ref{W1}) V
\\
 G \text{-}\mr{Op}(X^{D_{{\bf x}}\text{-} \mr{log}}; \chi (\check{\pmb{\lambda}})) @>> (\ref{Ee900}) > G \text{-}\mr{Op}(X^{D_{({\bf x}, {\bf z})}\text{-} \mr{log}}; (\chi (\check{\pmb{\lambda}}), \chi (-\check{\rho})^{m})).
\end{CD}
\end{align}
 If we  write
    \begin{align} \label{W21}
   G \text{-}\mr{MOp}(X^{D_{{\bf x}}\text{-} \mr{log}}; \check{\pmb{\lambda}})_{\mr{triv}, + {\bf w}} &:=  \coprod_{{\bf z} \in C^m ({\bf x})}    G \text{-}\mr{MOp}(X^{D_{{\bf x}}\text{-} \mr{log}};  \check{\pmb{\lambda}})_{\mr{triv}, +({\bf z}; {\bf w})},
  \end{align}
  then the upper horizontal arrow in  (\ref{W2})  in the case of  each ${\bf z} \in C^m ({\bf x})$ gives 
   an inclusion
  \begin{align} \label{W22}
  G \text{-}\mr{MOp}(X^{D_{{\bf x}}\text{-} \mr{log}}; \check{\pmb{\lambda}})_{\mr{triv}, + {\bf w}} \migiincl G \text{-}\mr{MOp}(X^{D_{{\bf x}}\text{-} \mr{log}}; \check{\pmb{\lambda}})_{\mr{gen}, + {\bf w}}.
  \end{align}
We shall consider 
$ G \text{-}\mr{MOp}(X^{D_{{\bf x}}\text{-} \mr{log}}; \check{\pmb{\lambda}})_{\mr{triv}, + {\bf w}}$ as a subset of $G \text{-}\mr{MOp}(X^{D_{{\bf x}}\text{-} \mr{log}}; \check{\pmb{\lambda}})_{\mr{gen}, + {\bf w}}$ via this injection.

\vspace{10mm}
\section{Dormant Miura opers and the Bethe ansatz equations}\SSP

In this section, 
we formulate and prove Theorem \ref{ThgA}  (cf. Theorem \ref{TT011}, Proposition  \ref{pgT011}).

\LSP
\subsection{Dormant Miura opers} \label{SS139}

The main result of this subsection 
 shows (cf. Proposition \ref{P015}) that each  dormant generic Miura oper specifies an element of  the set
$G \text{-}\mr{MOp}(X^{D_{{\bf x}}\text{-} \mr{log}}; \check{\pmb{\lambda}})_{\mr{triv}, + ({\bf z}; {\bf w})}$
 defined in the previous section (cf.  (\ref{Ee805})).

To begin with, recall 
from ~\cite[Definition 3.8.1]{Wak7}  
 that a Miura $G$-oper $\widehat{\mcE}^\spadesuit := (\mcE, \nabla, \mcE_B, \mcE'_B)$  is called   {\bf dormant} if $\nabla$ has vanishing $p$-curvature (cf., e.g., ~\cite[\S\,5.0]{K} or ~\cite[\S\,1.6]{Wak7}  for the  definition of $p$-curvature).
Then,   the following assertion holds.

\SSP
\bpr\label{P315} 
 Let  $\check{\lambda}$ be an element of $\mft$ such that there exists a dormant generic Miura $G$-oper   $\widehat{\mcE}^\spadesuit := (\mcE, \nabla,\mcE_B, \mcE'_B)$  on $\mbD^\mr{log}$ of exponent $\check{\lambda}$.
 Then,  $\check{\lambda}$ lies in $\mft_\mr{reg}^F$. 
 \epr
\begin{proof} 
Denote by $\mu \in \mfg$ the monodromy operator  (in the sense of ~\cite[Definition 1.46]{Wak5}) of $\nabla$ at the closed point $x_0$.
Also, denote by $\mr{ad} : \mfg \migi \mr{End}(\mfg)$ the adjoint representation of $\mfg$, which is injective and compatible with the respective natural restricted structures, i.e.,  $p$-power operations.
If we fix   a Jordan decomposition $\mu = \mu_s +\mu_n$ with $\mu_s$ semisimple  and $\mu_n$ nilpotent,
then there exists  an isomorphism $\alpha : \mr{End}(\mfg) \isom \mfg \mfl_{\mr{dim}(\mfg)}$ of restricted Lie algebras which sends
$\alpha (\mr{ad}(\mu))  \left(= \alpha (\mr{ad}(\mu_s)) + \alpha (\mr{ad}(\mu_n)) \right)$
to a Jordan normal form in $\mfg \mfl_{\mr{dim}(\mfg)}$.
In particular, $\alpha (\mr{ad}(\mu_s))$ is diagonal and every entry of $\alpha (\mr{ad}(\mu_n))$ except the superdiagonal  is $0$.
Let us observe the following sequence of equalities:
\begin{align} \label{Ww100}
\alpha (\mr{ad} (\mu_s)) + \alpha (\mr{ad}(\mu_n)) & = \alpha (\mr{ad} (\mu_s + \mu_n)) = \alpha (\mr{ad}((\mu_s +\mu_n)^{[p]})) \\
& = \alpha (\mr{ad}(\mu_s + \mu_n))^p = (\alpha (\mr{ad}(\mu_s))+ \alpha (\mr{ad}(\mu_n)))^p,\notag 
\end{align}
where the second equality follows from the assumption that $\nabla$ has vanishing $p$-curvature and $(-)^{[p]}$ denotes the $p$-power operation on $\mfg$ (cf. ~\cite[\S\,3.2.3 and \S\,3.4.2]{Wak5}).
By an explicit computation of 
$(\alpha (\mr{ad}(\mu_s))+ \alpha (\mr{ad}(\mu_n)))^p$, (\ref{Ww100}) implies that $\alpha (\mr{ad}(\mu_n)) =0$ (hence $\mu_n=0$), namely, $\mu$ is conjugate to some $v \in \mft$.
On the other hand,  $\mu$ is, by definition,  conjugate to $p_{-1} + \check{\lambda}$.
It follows that 
$p_{-1} + \check{\lambda}$ is conjugate to $v$, and hence, 
$\check{\lambda} = w (v)$ for some $w \in W$.
Since $p_{-1}+ \check{\lambda}$, as well as $v$,  is regular,
 $\check{\lambda}  \left(= w (v) \right)$ turns out to be regular.
Moreover, by the equality $\mu_n =0$, (\ref{Ww100})  reads the equality $\alpha (\mr{ad}(\mu)) = \alpha (\mr{ad} (\mu))^p$.
This  implies the equality  $v^{[p]}= v$, or equivalently,  $v \in \mft^F_\mr{reg}$.
This  completes the proof of the assertion.
\end{proof}
\SSP

Let $X$ be a connected  proper smooth curve over $k$, ${\bf x} := (x_i)_{i=1}^r$ ($r \geq 1$)  an ordered collection of distinct closed points of $X$. 
Also,  let $\check{\pmb{\lambda}} := (\check{\lambda}_i)_{i=1}^r \in (\mft_\mr{reg})^{r}$,  ${\bf w} \in W^{m}$.
Given each ${\bf z} \in C^m ({\bf x})$, we  obtain the subset
\begin{align}
   G \text{-}\mr{MOp}(X^{D_{({\bf x}, {\bf z})}\text{-}\mr{log}}; (\check{\pmb{\lambda}}, {\bf w}(-\check{\rho})))_{\mr{gen}}^{^\mr{Zzz...}}
\end{align}
 of  
$ G \text{-}\mr{MOp}(X^{D_{({\bf x}, {\bf z})}\text{-}\mr{log}}; (\check{\pmb{\lambda}}, {\bf w}(-\check{\rho})))_{\mr{gen}}$
 consisting of dormant generic Miura $G$-opers.
 According to  the above proposition, this subset is empty unless $\check{\pmb{\lambda}} \in (\mft_\mr{reg}^F)^{r}$.

Also, we write
\begin{align} \label{Ee600}
 G \text{-}\mr{MOp}(X^{D_{{\bf x}}\text{-}\mr{log}}; \check{\pmb{\lambda}})_{\mr{gen}, +{\bf w}}^{^\mr{Zzz...}} 
 := \coprod_{{\bf z}\in C^m ({\bf x})} 
  G \text{-}\mr{MOp}(X^{D_{({\bf x}, {\bf z})}\text{-}\mr{log}}; (\check{\pmb{\lambda}}, {\bf w}(-\check{\rho})))_{\mr{gen}}^{^\mr{Zzz...}}.
\end{align}
For an element $v \in \mfg$ and a positive integer $s$, we shall set
\begin{align} \label{Ww552}
v^s := (v,v, \cdots, v) \in \mfg^s.
\end{align}
Then, the following assertion holds.

\SSP
\bpr \label{P015}
Let us keep the above notation.
\begin{itemize}
\item[(i)]
The subset $G \text{-}\mr{MOp}(X^{D_{({\bf x}, {\bf z})}\text{-}\mr{log}}; (\check{\pmb{\lambda}}, {\bf w}(-\check{\rho})))_{\mr{gen}}^{^\mr{Zzz...}}$ of $G \text{-}\mr{MOp}(X^{D_{({\bf x}, {\bf z})}\text{-}\mr{log}}; (\check{\pmb{\lambda}}, {\bf w}(-\check{\rho})))_{\mr{gen}}$  
 (for each  ${\bf z} \in C^m ({\bf x})$)
 is contained in the subset 
$G \text{-}\mr{MOp}(X^{D_{{\bf x}}\text{-} \mr{log}}; \check{\pmb{\lambda}})_{\mr{triv}, + ({\bf z}; {\bf w})}$.
In particular, we have an inclusion
\begin{align} \label{W250}
 G \text{-}\mr{MOp}(X^{D_{{\bf x}}\text{-}\mr{log}}; \check{\pmb{\lambda}})_{\mr{gen}, +{\bf w}}^{^\mr{Zzz...}} 
\migiincl 
G \text{-}\mr{MOp}(X^{D_{{\bf x}}\text{-}\mr{log}}; \check{\pmb{\lambda}})_{\mr{triv}, + {\bf w}}. 
\end{align}
\item[(ii)]
Suppose further that $\check{\pmb{\lambda}} = (-\check{\rho})^{r}$ and that  $X$ is the projective line $\mbP$ over $k$.
Then, the injection (\ref{W250}) becomes bijective.
\end{itemize}
 \epr
\begin{proof} 
First, we shall consider assertion (i).
Let $w \in W$ and  let $\widehat{\mcE}^\spadesuit := (\mcE, \nabla, \mcE_B, \mcE'_B)$ be  a dormant  generic Miura $G$-oper 
on $\mbD^\mr{log}$ of exponent $w (-\check{\rho})$.
Let us choose a uniformizer $t$ of $\widehat{\mcO}$, which induces  $\widehat{\mcO} \isom k[\![t]\!]$.
To complete the proof, it suffices to prove that 
$(\mcE, \nabla, \mcE_B)$ becomes a $G$-oper on $\mbD$ after a gauge transformation by some element of $B (k(\!(t)\!))$.
By the assumption $(*)_G$,
one may obtain the exponential map $\mr{exp} : \mfn \migi N $ given by  ~\cite[Proposition 1.31 and Remark 1.35]{Wak5}.
This map enables us to apply 
 an argument similar to  the proof of ~\cite[Proposition 9.2.1]{Fr} to our positive characteristic case.
Hence, 
after a gauge transformation by some element of $B (k(\!(t)\!))$,
$\nabla$ may be expressed as 
$\partial_t + p_{-1}+ {\bf v}(t) + \frac{v}{t}$ for some ${\bf v}(t) \in \mfb [\![t]\!]$ and $v \in \mfn$.
The mod $t$ reduction of the  $p$-curvature of $\nabla$  is given by $v^{[p]} - v$, which
 is equal to 
$0$  because of the dormancy condition on $\widehat{\mcE}^\spadesuit$.
But, 
 since 
  $v \in \mfn$, we have $v^{[p]}=0$,  which implies $v = 0$.
Therefore,  $(\mcE, \nabla, \mcE_B)$ forms a $G$-oper on $\mbD$, and this completes the proof of assertion (i).

Next, we shall  consider assertion (ii).
Let  $\widehat{\mcE}^\spadesuit  := (\mcE, \nabla, \mcE_B, \mcE'_B)$ be a Miura $G$-oper  classified by 
  $G \text{-} \mr{MOp} (\mbP^{D_{{\bf x}}\text{-} \mr{log}}; (-\check{\rho})^{r})_{\mr{triv}, + {\bf w}}$, i.e., by  $G \text{-} \mr{MOp} (\mbP^{D_{{\bf x}}\text{-} \mr{log}};  (-\check{\rho})^{r})_{\mr{triv}, + ({\bf z};  {\bf w})}$
   for some ${\bf z}:= (z_1, \cdots, z_m)\in C^m ({\bf x})$.
By applying the discussion in Remark \ref{R1019} to the restriction $\widehat{\mcE}^\spadesuit |_{\mbD_{x_i}}$ (for each $i=1, \cdots, r$), we can verify
 that
$\widehat{\mcE}^\spadesuit$ comes, via gauge transformation,  from a Miura $G$-oper on $\mbP^{D_{{\bf z}}\text{-} \mr{log}}$ of exponents ${\bf w}$.
Moreover,
since
the restriction of $\widehat{\mcE}^\spadesuit$   to $\mbD_{z_j}$ (for each $j = 1, \cdots, m$) belongs to  
$G \text{-} \mr{Op}(\mbD_{z_j})  \left(\subseteq G \text{-} \mr{Op} (\mbD_{z_j}^\mr{log}; \chi (-\check{\rho})) \right)$,
 this $G$-oper comes from a $G$-oper $\mcE^\spadesuit_0$ on $\mbP$.
By  ~\cite[Corollary 2.6.2 and Theorem 3.3.1, (iv)]{Wak6},
 $\mcE^\spadesuit_0$ turns out to be  dormant.
It follows that $\widehat{\mcE}^\spadesuit$ is dormant, and hence,
belongs to 
$ G \text{-} \mr{MOp} (\mbP^{D_{{\bf x}}\text{-} \mr{log}};  (-\check{\rho})^{r})^{^\mr{Zzz...}}_{\mr{gen}, + {\bf w}}$.
This completes the proof of  assertion (ii).
\end{proof}

\LSP
\subsection{Bethe ansatz equations modulo $p$} \label{SS131gg}

We here consider the case where $X$ is  taken to be  the projective line $\mbP$  over $k$.
 Denote by $x$ the natural coordinate of $\mbP$, i.e.,  $\mbP \setminus \{ \infty \}= \mr{Spec}(k[x])$.
 Fix integers $r$, $m$  with $r \geq 0$, $m \geq 1$.
Let  ${\bf x} := (x_1, \cdots, x_r, x_{r+1})$ be   an ordered collection of distinct closed points of $\mbP$ with $x_{r+1}=\infty$, and let 
 $\check{\pmb{\lambda}}^0 := (\check{\lambda}_1^0, \cdots, \check{\lambda}^0_r, \check{\lambda}_{r+1}^0) \in (\mft_\mr{reg})^{r+1}$, 
 $\check{\pmb{\lambda}}' := (\check{\lambda}'_1, \cdots, \check{\lambda}'_m) \in (\mft_\mr{reg})^{m}$.

\SSP
\bpr \label{P05} 
Let  ${\bf z} \in C^m ({\bf x})$.
Then, 
the set $G \text{-}\mr{MOp}(\mbP^{D_{({\bf x}, {\bf z})}\text{-}\mr{log}}; (\check{\pmb{\lambda}}^0, \check{\pmb{\lambda}}'))_\mr{gen}$ is nonempty if and only if the equality
\begin{align} \label{Ee907}
\sum_{i=1}^{r+1} (\check{\lambda}_{i}^0+\check{\rho}) +\sum_{j=1}^m (\check{\lambda}'_j +\check{\rho}) = 2\check{\rho}
\end{align}
holds.
Moreover, if  the  equality (\ref{Ee907}) holds, then 
$G \text{-}\mr{MOp}(\mbP^{D_{({\bf x}, {\bf z})}\text{-}\mr{log}}; (\check{\pmb{\lambda}}^0, \check{\pmb{\lambda}}'))_\mr{gen}$
 consists exactly of a single  element.
 \epr
\begin{proof} 
Let us set 
 $\mcE^\dagger_{T, \mr{log}} := \Omega_{\mbP^{D_{({\bf x}, {\bf z})}\text{-}\mr{log}}/k}^\times \times^{\mbG_m, \check{\rho}} T$ and $\mcE^\dagger_T :=\Omega_{\mbP/k}^\times \times^{\mbG_m, \check{\rho}} T$, where for each line bundle $\mcL$ we denote by $\mcL^\times$ the $\mbG_m$-bundle corresponding to $\mcL$.
Recall from ~\cite[Proposition 3.7.1]{Wak7} that 
there exists a canonical bijection
\begin{align} \label{Ee905}
G \text{-}\mr{Conn}(\mbP^{D_{({\bf x}, {\bf z})}\text{-}\mr{log}}; (\check{\pmb{\lambda}}^0, \check{\pmb{\lambda}}')) 
\isom
 G \text{-}\mr{MOp}(\mbP^{D_{({\bf x}, {\bf z})}\text{-}\mr{log}}; (\check{\pmb{\lambda}}^0, \check{\pmb{\lambda}}'))_\mr{gen}, 
\end{align}
where the left-hand side denotes the  set of $D_{({\bf x}, {\bf z})}$-log connections on $\mcE^\dagger_{T, \mr{log}}$  whose  monodromy operators are given by $(\check{\pmb{\lambda}}^0, \check{\pmb{\lambda}}')$.

Suppose that there exists a $D_{({\bf x}, {\bf z})}$-log connection  $\overline{\nabla}$ classified by 
$G \text{-}\mr{Conn}(\mbP^{D_{({\bf x}, {\bf z})}\text{-}\mr{log}}; (\check{\pmb{\lambda}}^0, \check{\pmb{\lambda}}'))$.
One may find  a unique $D_{({\bf x}, {\bf z})}$-log  connection $\overline{\nabla}'$ on
$\mcE^\dagger_T$
whose restriction to the curve $\mbP^\circledcirc  := \mbP \setminus \{ x_1,  \cdots, x_{r+1}, z_1, \cdots, z_m \}$ 
coincides with  $\overline{\nabla} |_{\mbP^\circledcirc}$.
The monodromy operators of $\overline{\nabla}'$ are  $(\check{\lambda}_1^0+\check{\rho}, \cdots,  \check{\lambda}^0_{r+1}+\check{\rho}, \check{\lambda}'_1+\check{\rho}, \cdots,  \check{\lambda}'_m+\check{\rho})$, so  $\overline{\nabla}'$ must  be
  expressed as 
 \begin{align} \label{Ee901}
 \partial_x 
  + \sum_{i=1}^r \frac{\check{\lambda}_i^0+ \check{\rho}}{x - x_i} + \sum_{j=1}^m \frac{\check{\lambda}'_j + \check{\rho}}{x -z_j},
\end{align}
where $\partial_x = \frac{d}{d x}$, under 
the trivialization of the $T$-bundle $\mcE_T^\dagger |_{\mbP\setminus \{ \infty \}}$
given by $\Omega_{\mbP/k} |_{\mbP \setminus \{ \infty\}} \isom \mcO_{\mbP \setminus \{ \infty\}}; d x \mapsto 1$ (cf. \cite[\S\,3.1,   (3.1)]{Fre2}).
 In particular, $G \text{-}\mr{Conn}(\mbP^{D_{({\bf x}, {\bf z})}\text{-}\mr{log}}; (\check{\pmb{\lambda}}^0, \check{\pmb{\lambda}}'))$ consists of a single  connection.
Moreover, 
according to 
the discussion in {\it loc.\,cit.},
we have the equality
\begin{align}
2 \check{\rho} -\sum_{i=1}^r (\check{\lambda}_i^0+\check{\rho}) -\sum_{j=1}^m (\check{\lambda}'_j +\check{\rho}) = \check{\lambda}_{r+1}^0 +\check{\rho},
\end{align}
which is equivalent to  the  equality (\ref{Ee907}).

Conversely, if   the equality (\ref{Ee907}) holds, 
then 
we obtain  a unique $D_{({\bf x}, {\bf z})}$-log  connection $\overline{\nabla}'$ on 
$\mcE^\dagger_T$
 of the  form  
(\ref{Ee901}).
Moreover, there exists a unique $D_{({\bf x}, {\bf z})}$-log connection $\overline{\nabla}$ on $\mcE_{T, \mr{log}}^\dagger$ whose restriction to 
$\mbP^\circledcirc$
coincides  with $\overline{\nabla}' |_{\mbP^\circledcirc}$.
Since $\overline{\nabla}$
  belongs to 
   $G \text{-}\mr{Conn}(\mbP^{D_{({\bf x}, {\bf z})}\text{-}\mr{log}}; (\check{\pmb{\lambda}}^0, \check{\pmb{\lambda}}'))$, 
the bijection  (\ref{Ee905}) implies  that   $G \text{-}\mr{MOp}(\mbP^{D_{({\bf x}, {\bf z})}\text{-}\mr{log}}; (\check{\pmb{\lambda}}^0, \check{\pmb{\lambda}}'))_\mr{gen}$ is nonempty.
This completes the proof of the assertion.
\end{proof}
\SSP

Let  $\pmb{\alpha} := (\alpha_1, \cdots, \alpha_m) \in \Gamma^{m}$.
For each $j =1, \cdots, m$, we shall  write $\alpha^W_j \in W$ for the simple reflection corresponding to $\alpha_j$.
The coroot $\check{\alpha}_j$ corresponding to $\alpha_j$ satisfies the equality  $\check{\alpha}_j =  \alpha^W_j (-\check{\rho})+\check{\rho}$.
Also, we shall write $\pmb{\alpha}^W := (\alpha^W_1, \cdots, \alpha^W_m)$.
Next, let us take an element $\check{\pmb{\lambda}} := (\lambda_1, \cdots, \lambda_{r+1})$ of  $(\mft^F_\mr{dom})^{r+1}$ satisfying the equality 
\begin{align} \label{W100}
-\sum_{i=1}^{r+1} \check{\lambda}_i + \sum_{j=1}^m \check{\alpha}_j = 2 \check{\rho}
\end{align}
(i.e., the equality (\ref{Ee907}) in the case where the data $(\check{\pmb{\lambda}}^0, \check{\pmb{\lambda}}')$ is replaced with  $(-\check{\pmb{\lambda}}-\check{\rho}, \pmb{\alpha}^W (-\check{\rho}))$).
Note that (since $\check{\pmb{\lambda}} \in (\mft^F_\mr{dom})^{r+1}$) the element 
$-\check{\pmb{\lambda}}-\check{\rho} := (-\check{\lambda}_1 -\check{\rho}, \cdots, -\check{\lambda}_{r+1} -\check{\rho})$ belongs to $(\mft^F_\mr{reg})^{r+1}$.
The assignment from each ${\bf z} \in C^m ({\bf x})$
to a unique generic Miura $G$-oper on $\mbP^{D_{({\bf x}, {\bf z})}\text{-}\mr{log}}$ of exponents $(-\check{\pmb{\lambda}}-\check{\rho}, \pmb{\alpha}^W (-\check{\rho}))$ (cf. Proposition \ref{P05}) defines a bijection
 \begin{align} \label{Ee045}
 C^m ({\bf x}) \isom  G \text{-} \mr{MOp} (\mbP^{D_{{\bf x}}\text{-}\mr{log}}; -\check{\pmb{\lambda}}-\check{\rho})_{\mr{gen}, + \pmb{\alpha}^W},
 \end{align}
 where the right-hand side contains $G \text{-} \mr{MOp} (\mbP^{D_{{\bf x}}\text{-}\mr{log}}; -\check{\pmb{\lambda}}-\check{\rho})_{\mr{triv}, + \pmb{\alpha}^W}$ (cf. (\ref{W22})).
 
   Now, denote by
\begin{align} \label{Ef22}
G \text{-} \mr{BA}_{\check{\pmb{\lambda}}, \pmb{\alpha}}  \left( \subseteq C^m ({\bf x}) \right)
\end{align}
the set of elements ${\bf z} := (z_1, \cdots, z_m)$ of $C^m ({\bf x})$ satisfying  
the system of equations
\begin{align} \label{Ee070}
\sum_{i=1}^r \frac{\langle \alpha_j, \check{\lambda}_i\rangle}{z_j - x_i} - \sum_{s \neq j} \frac{\langle \alpha_j, \check{\alpha}_s\rangle}{z_j -z_s} =0 \hspace{5mm} \left(j=1, \cdots, m\right),
\end{align}
i.e., 
the {\bf Bethe ansatz equations} (cf. Introduction).

\SSP
\begin{rema}\label{Rr319}
We shall consider 
the case where  $G = \mr{PGL}_2$ and $p \geq 3$  (i.e.,   the condition   $(*)_{\mr{PGL}_2}$ is satisfied).
Then, we have 
 $\check{\rho} = \begin{pmatrix} \frac{1}{2} & 0 \\ 0 & -\frac{1}{2}\end{pmatrix}$,  $\Gamma = \{ \alpha \}$, and  $\check{\alpha} = 2 \check{\rho} = \begin{pmatrix} 1 & 0 \\ 0 & -1\end{pmatrix}$.
For each $i \in \{1, \cdots r \}$, denote by $q_i$ the element of $k$ with
 $\check{\lambda}_i = q_i \cdot \check{\rho}$.
Then,  (\ref{Ee070}) reads the system of  equations 
\begin{align}
\label{We070}
\sum_{i=1}^r \frac{q_i}{z_j - x_i} - \sum_{s \neq j} \frac{2}{z_j -z_s} =0 \hspace{5mm} \left(j=1, \cdots, m\right)
\end{align}
of values in $k$.

In what follows, we suppose further that $r =0$.
In particular, (\ref{We070})  is equivalent to 
\begin{align} \label{Ee807}
\sum_{s \neq j}\frac{1}{z_j -z_s} =0\hspace{5mm} \left(j=1, \cdots, m\right).
\end{align}
 Write $f(x)$ for the polynomial defined as
 \begin{align} \label{Ee8001}
 f (x) := \prod_{j=1}^{m} (x -z_j) \in k[x].
 \end{align}
Since $z_1, \cdots, z_{m}$ are distinct, we have $\mr{gcd}(f (x), f'(x)) =1$, or equivalently, $f'(z_j) \neq 0$ for every  $j \in \{1, \cdots, m \}$.
Also,  for each  $j \in \{1, \cdots, m \}$, the following sequence of equalities holds:
\begin{align}
\frac{f''(z_j)}{f'(z_j)} & = \frac{\sum_{s'=1}^{m} \sum_{s'' \neq s'} \prod_{s'''\neq s', s''} (x - z_{s'''}) |_{x = z_j}}{\sum_{s'=1}^{m} \prod_{s'' \neq s'} (x - z_{s''}) |_{x = z_j}}  \\
& = \frac{2 \cdot \sum_{s'\neq s}\prod_{s'' \neq j, s'} (z_j -z_{s''})}{\prod_{j \neq s}(z_j - z_s)} \notag \\
& = 2 \cdot \sum_{s\neq j} \frac{1}{z_j-z_s}. \notag
\end{align}
On the other hand, (since  $f (z_j) =0$ for any $j \in \{1, \cdots, m \}$)  the equality  $\frac{f''(z_j)}{f'(z_j)} = 0$
holds  for any $j \in \{1, \cdots, m \}$ if and only if
$f''(x) =0$ as an element of $k[x]$.
 Thus,  we obtain the following equivalence of conditions:
 \begin{align} \label{Ee801}
 \sum_{s\neq j} \frac{1}{z_j-z_s} =0 \hspace{3mm} \left(\text{for every} \ j \in \{1, \cdots, m \}\right) \ \ 
 \Longleftrightarrow  \ \  f''(x) =0. 
 \end{align}
 
 For instance, if $(a,b,c)$ is  an element of $k^3$ with $ab-c \neq 0$ and   $z_{1}, \cdots, z_{p+1}$ are  the roots of the polynomial $f(x) := x^{p+1} + a x^p + b x + c$,
 then the equalities  $\sum_{s\neq j} \frac{1}{z_j - z_s}=0$ ($j =1, \cdots, p+1$) are satisfied.

More concretely,  if we   consider the case where $p=3$ and  $m=4$,  then 
\begin{align}
(z_1, z_2, z_3, z_4) = \left(1 + 2 \sqrt{1 + 2 \sqrt{2}}, 1 - 2 \sqrt{1 + 2 \sqrt{2}},  1 +  \sqrt{2 +  \sqrt{2}}, 1 -  \sqrt{2 +  \sqrt{2}}\right)
\end{align}
specifies a solution to the system of equations (\ref{Ee807}).
\end{rema}

\LSP
\subsection{Comparison I}
  \label{SS137}

In this subsection, we prove the following Theorem \ref{TT011}, being a part of Theorem A.
The corresponding  assertion 
in the complex (i.e., original) case
was proved in   ~\cite[Theorem 3.2]{Fre2}.
Also, Lemma \ref{L02}, which is  described later and used in the  proof of Theorem \ref{TT011},
corresponds to ~\cite[Lemma 2.10]{Fre2}.
Unlike the proofs  in {\it loc.\,cit}.,  we will not make  any  analytic argument 
in order to  include the case of positive characteristic.
    
\SSP
\bt \label{TT011} 
 Let us keep the notation in the discussion  following Proposition \ref{P05}, and suppose that 
 the equality (\ref{W100}) holds.
Then, the bijection 
 (\ref{Ee045})
  restricts to a bijection
 \begin{align} \label{Ee55}
 G \text{-} \mr{BA}_{\check{\pmb{\lambda}}, \pmb{\alpha}} 
 \isom
 G \text{-} \mr{MOp} (\mbP^{D_{{\bf x}}\text{-}\mr{log}}; -\check{\pmb{\lambda}} -\check{\rho})_{\mr{triv}, + \pmb{\alpha}^W}.
 \end{align}
 If, moreover,  $\check{\pmb{\lambda}} = 0^{r+1}$, then the bijection (\ref{Ee55})  becomes the following bijection:
 \begin{align} \label{Ee990}
 G \text{-} \mr{BA}_{0^{r+1}, \pmb{\alpha}} 
 \isom
 G \text{-} \mr{MOp} (\mbP^{D_{{\bf x}}\text{-}\mr{log}}; (-\check{\rho})^{r+1})_{\mr{gen}, + \pmb{\alpha}^W}^{^\mr{Zzz...}}.
 \end{align}
  \et
\begin{proof}
Since the latter assertion follows directly from the former assertion and Proposition \ref{P015}, (ii), it suffices to prove the former assertion.
Let ${\bf z} :=(z_1, \cdots, z_m)$ be an element of $C^m ({\bf x})$ and denote by 
$\widehat{\mcE}^\spadesuit := (\mcE, \nabla, \mcE_B, \mcE'_B)$  the Miura $G$-oper 
corresponding to the  unique element of 
$G \text{-} \mr{MOp} (\mbP^{D_{({\bf x}, {\bf z})}\text{-}\mr{log}}; -\check{\pmb{\lambda}} -\check{\rho}, \pmb{\alpha}^W (-\check{\rho}))_{\mr{gen}}$.
By considering the discussion in the proof of Proposition \ref{P05},
we see that $\nabla |_{\mbP \setminus \{ \infty \}}$ may be expressed as 
\begin{align}
\partial_x + p_{-1} - \sum_{i=1}^r \frac{\check{\lambda}_i}{x -x_i} + \sum_{j=1}^m \frac{\check{\alpha}_j}{x - z_j}
\end{align} 
under  a suitable trivialization  of $\mcE |_{\mbP \setminus \{ \infty \}}$.
 Then,  
 the restriction of $\nabla$ to $\mbD_{z_j}$ (for each  $j \in \{1, \cdots, m \}$) may be expressed, after choosing a uniformizer $t$ at $z_j$, as 
 \begin{align}
 \partial_t + p_{-1} + \frac{\check{\alpha}_j}{t} + {\bf u}_j (t),
 \end{align}
where $\partial_t := \frac{d}{dt}$ and ${\bf u}_j (t)$  denotes   a certain element of $\mft [\![t]\!]$.
 One verifies that
 the equality $\langle \alpha_j, {\bf u}_j (0) \rangle =0$ holds
if and only if ${\bf z}$ satisfies  the $j$-th equation in (\ref{Ee070}).
Thus, the assertion follows from the following lemma. 
\end{proof}

\SSP
\ble \label{L02}
Let $\alpha \in \Gamma$, and denote by $\alpha^W$ the simple reflection corresponding to $\alpha$.
Choose a uniformizer $t$ of $\widehat{\mcO}$.
Also, let us take a $G$-oper  $\mcE^\spadesuit$ on $\mbD^\mr{log}$ of the form  $(\mbD \times G, \nabla, \mbD \times B)$, where
$\nabla = \partial_t + p_{-1} + \frac{\check{\alpha}}{t} +
   {\bf u}(t)$ (under the identification $k[\![t]\!] \isom \widehat{\mcO}$ determined  by $t$) for some ${\bf u}(t) \in \mft [\![t]\!]$.
Then, 
$\langle \alpha, {\bf u}(0) \rangle =0$ if and only if $\mcE^\spadesuit$ belongs to $G \text{-} \mr{Op}(\mbD)$ (cf. (\ref{Ee8801})).
 \ele
\begin{proof}
To begin with, let us make the following observation.
Denote by $e_\alpha$ the unique generator   of $\mfg^\alpha$ such
that $\{ f_\alpha, 2 \check{\alpha}, e_\alpha \}$ forms an $\mfs \mfl_2$-triple.
Also, denote by $\mr{exp} : \mfn \migi N$ the exponential map asserted in  ~\cite[Proposition 1.31]{Wak5}.
It follows from   ~\cite[Corollary 1.34]{Wak5} that after  the gauge transformation by $\mr{exp}(-\frac{1}{t} \cdot e_\alpha)^{-1}$, $\nabla$ becomes the connection $\nabla'$ of the  form
\begin{align} \label{W200}
\nabla' := \partial_t +  \frac{d }{d t} \left(\frac{1}{t} \cdot e_\alpha \right) + \sum_{s=0}^\infty \frac{1}{s!} \cdot \mr{ad} \left(- \frac{1}{t} \cdot e_\alpha \right)^s \left(p_{-1}+ \frac{\check{\alpha}}{t} + {\bf u}(t) \right).
\end{align}
Observe that  $\frac{d}{d t} (\frac{1}{t} \cdot e_\alpha) = -\frac{1}{t^2}\cdot e_\alpha$ and
\begin{align}
\mr{ad} \left(- \frac{1}{t} \cdot e_\alpha \right) \left(p_{-1}+ \frac{\check{\alpha}}{t} + {\bf u}(t) \right) 
& = 
-\frac{1}{t} \cdot \check{\alpha} +  \frac{2}{t^2} \cdot e_\alpha - \frac{\langle \alpha, {\bf u}(t) \rangle}{t} \cdot e_\alpha,  \\ 
\mr{ad} \left(- \frac{1}{t} \cdot e_\alpha \right)^2 \left(p_{-1}+ \frac{\check{\alpha}}{t} + {\bf u}(t) \right) & =
-\frac{2}{t^2} \cdot e_\alpha, \notag \\ 
\mr{ad} \left(- \frac{1}{t} \cdot e_\alpha \right)^l \left(p_{-1}+ \frac{\check{\alpha}}{t} + {\bf u}(t) \right) 
& =  0 \notag
\notag
\end{align}
($l =3,4, \cdots$).
Thus, we have the equality 
\begin{align}
\nabla' = \partial_t + p_{-1} + {\bf u}(t)- \frac{\langle \alpha, {\bf u}(t) \rangle}{t} \cdot e_\alpha.
\end{align}
Hence, if  $\langle \alpha, {\bf u}(0) \rangle =0$, then 
 the triple $(\mbD \times G, \nabla', \mbD \times B)$, as well as $\mcE^\spadesuit$,  belongs to $G \text{-} \mr{Op}(\mbD)$.
This completes
the ``only if'' part of the assertion.

Next, we shall suppose that 
$\mcE^\spadesuit$ belongs to $G\text{-}\mr{Op}(\mbD)$.
There exist  elements $b \in B (k(\!(t)\!))$ and  ${\bf u}_b(t) \in \mfb [\![t]\!]$  such that $\nabla'$ becomes the connection $\nabla'_b = \partial_t + p_{-1} + {\bf u}_b (t)$ after the  gauge transformation by $b$.
By an explicit computation of this gauge transformation, 
 $b$ turns out to be contained in   $N (k(\!(t)\!))$, i.e., $b = \mr{exp}({\bf v}(t))^{-1}$ for some  ${\bf v} (t) \in \mfn (\!(t)\!)$.
That is to say, we have an equality
\begin{align} \label{Ee22}
& \ \partial_t + p_{-1} + {\bf u}_b (t)   \left(= \nabla'_b \right)  \\
= & \  \partial_t + \frac{d}{dt} (-{\bf v}(t)) + \sum_{s=0}^\infty \frac{1}{s!} \cdot \mr{ad} ({\bf v}(t))^s \left(p_{-1} + {\bf u}(t) - \frac{\langle \alpha, {\bf u}(t) \rangle}{t} \cdot e_\alpha \right). \notag
\end{align}
According to the principal  gradation on $\mfg$,  
the elements ${\bf v}(t)$ and  ${\bf u}_b (t)$ can be decomposed as ${\bf v}(t) = \sum_{i=1}^{\mr{rk}(\mfg)} {\bf v}_i (t)$ and ${\bf u}_b (t) = \sum_{i=1}^{\mr{rk}(\mfg)} {\bf u}_{b, i} (t)$, respectively, where ${\bf u}_{b, i} (t)  \in \mfg_i [\![t]\!]$ and  ${\bf v}_i (t)\in \mfg_i (\!(t)\!)$.
Given a pair of positive integers $(i, j)$ with  $j \leq i \leq   \mr{rk}(\mfg)$, 
we denote by $F_j^i$ the degree $i$ graded $k$-linear endomorphism  of $\mfg (\!(t)\!)$ defined as
\begin{align}
F_j^i := \sum_{s=j}^i \frac{1}{s!} \cdot \sum_{\substack{(l_1, \cdots, l_s) \in \mbZ_{>0}^s, \\ l_1 + \cdots + l_s= i}} \mr{ad} ({\bf v}_{l_1}(t)) \circ \cdots \circ \mr{ad}({\bf v}_{l_s}(t)) : \mfg (\!(t)\!) \migi \mfg (\!(t)\!).
\end{align}
Then, the equality (\ref{Ee22}) implies the equality 
\begin{align} \label{Ee24}
{\bf u}_{b, i} (t) =  - \frac{d}{dt} ({\bf v}_i (t)) -  \frac{\langle \alpha, {\bf u}(t) \rangle}{t} \cdot  F_1^{i-1}\left( e_\alpha \right)+ F_1^i ({\bf u}(t)) + F_2^{i+1}(p_{-1}) + [{\bf v}_{i+1}(t), p_{-1}]
\end{align}
 for each $i=1, \cdots, \mr{rk}(\mfg)-1$ and the equality
 \begin{align} \label{Ee23}
 {\bf u}_{b, \mr{rk}(\mfg)} (t) =  - \frac{d}{dt} ({\bf v}_{ \mr{rk}(\mfg)} (t)) -  \frac{\langle \alpha, {\bf u}(t) \rangle}{t} \cdot  F_1^{ \mr{rk}(\mfg)-1}\left( e_\alpha \right)+ F_1^{ \mr{rk}(\mfg)} ({\bf u}(t)). 
 \end{align}
Let us assume  the condition ${\bf v}(t) \notin \mfn [\![t]\!]$, which implies that the positive integer 
\begin{align}
i_\mr{min} :=\mr{min} \left\{ i \in \mbZ_{>0} \, | \, {\bf  v}_i (t) \notin \mfg_i [\![t]\!] \right\}
\end{align}
 is well-defined.
Since $[-, p_{-1}]  \left(= - \mr{ad} (p_{-1}) \right) : \mfg_{i+1} \migi \mfg_i$ (where $i \geq 0$) is injective, 
the equality (\ref{Ee23}) implies, by induction on $i$, that
the pole order of ${\bf v}_{i_\mr{min} +j}(t)$ (for $j =1,2, \cdots $) coincides with  $1 + j$.
In particular, the pole order of ${\bf v}_{\mr{rk}(\mfg)}(t)$ is $1 + \mr{rk}(\mfg) - i_{\mr{min}}$.
But, by comparing the respective pole orders of the both sides of (\ref{Ee23}),
we have $\frac{d}{dt} ({\bf v}_{\mr{rk}(\mfg)}(t)) = 0$, which  contradicts   the inequality $1 + \mr{rk}(\mfg) - i_{\mr{min}} < p$ induced by   
the assumption $(*)_G$.
Therefore,  
${\bf v}(t)$ must be contained in $\mfn [\![t]\!]$, and hence, 
 $\nabla'$ forms a connection on $\mbD \times G$.
 It follows that 
 $\langle \alpha, {\bf u}(t) \rangle = 0$, as desired.
 This completes the proof of the lemma.
\end{proof}

\SSP
\begin{rema}\label{Rr119}
In ~\cite{Var},
A. Varchenko  studied the Bethe ansatz equations (for $\mfg = \mfs \mfl_2$) modulo $p$ and showed that
the Bethe vector corresponding to its  solution is an eigenvector of the Gaudin Hamiltonian.
By  this  result together with Theorem  \ref{TT011},
we can construct  eigenvectors by means of dormant generic Miura $\mr{PGL}_2$-opers.
On the other hand, as explained in  \S\,\ref{SS131} later,
dormant generic Miura $\mr{PGL}_2$-opers correspond bijectively to  (pre-)Tango structures.
In this way,  a certain type of (pre-)Tango structures provides eigenvectors of the Gaudin Hamiltonian.
 \end{rema}

\LSP
\subsection{Comparison II} \label{SS731}

In this subsection, we shall characterize, via the  bijection (\ref{Ee55}) obtained in the previous subsection,  the subsets of  $G \text{-}\mr{MOp}(\mbP^{D_{{\bf x}}\text{-}\mr{log}}; - \check{\pmb{\lambda}}-\check{\rho})_{\mr{triv}, + \pmb{\alpha}^W}$ consisting of Miura $G$-opers with vanishing $p$-curvature and   nilpotent $p$-curvature, respectively.
Let us keep the notation in \S\,\ref{SS131gg},
 and suppose that the equality (\ref{W100}) holds.

First, consider the case of Miura $G$-opers vanishing $p$-curvature, i.e., dormant Miura $G$-opers.
Given  an element $A$ of $\mfg (k (x))  \left(= \mfg \otimes_k k(x) \right)$,
we shall define elements $A (l)$ ($l =1,2, \cdots$) of $\mr{End}_k(\mfg (k (x)))$ determined inductively  by the following rule:
\begin{align} \label{Ee109}
A (1) := \mr{ad}(A), \hspace{5mm} A(l+1) := \partial_x (A (l)) + \mr{ad}(A) \circ   A(l) \ (l \geq 1).
\end{align}

Let $\nabla$ be the connection on the trivial $G$-bundle defined on  a suitable open subscheme  of $\mbP \setminus \{ \infty \}  \left(= \mr{Spec}(k[x]) \right)$ given by $\nabla = \partial_x + A$.
The $p$-curvature of $\nabla$ may be expressed, by definition,  as  $A (p)$ under the trivialization $\Omega_{\mbP \setminus \{ \infty \}} \isom \mcO_{\mbP \setminus \{ \infty \}}$ given by $dx \mapsto 1$.
In particular, $\nabla$ has vanishing $p$-curvature if and only if $A (p) =0$.

Now, let us write
\begin{align}
G \text{-} \mr{BA}_{\check{\pmb{\lambda}}, \pmb{\alpha}}^\circledast
\end{align}
for the subset of $G\text{-} \mr{BA}_{\check{\pmb{\lambda}}, \pmb{\alpha}}$ consisting of 
elements  ${\bf z} := (z_1, \cdots, z_m)$ with $A^{{\bf z}}(p) =0$, where 
\begin{align} \label{Ee110}
A^{{\bf z}} := p_{-1} -\sum_{i=1}^r \frac{\check{\lambda}_i}{x-x_i} + \sum_{j=1}^m \frac{\check{\alpha}_i}{x -z_j} \in \mfg (k(x)).
\end{align}
Then, the above  discussion implies the following proposition.

\SSP
\bpr \label{pggT1} 
 The bijection (\ref{Ee55}) restricts to a bijection 
 \begin{align}
G \text{-} \mr{BA}_{\check{\pmb{\lambda}}, \pmb{\alpha}}^\circledast \isom 
G \text{-}\mr{MOp}(\mbP^{D_{{\bf x}}\text{-}\mr{log}}; - \check{\pmb{\lambda}}-\check{\rho})_{\mr{gen}, + \pmb{\alpha}^W}^{^\mr{Zzz...}}.
 \end{align}
  \epr
\SSP

Next, we shall consider the case of Miura $G$-opers with nilpotent $p$-curvature.
Denote by $[0]_\mfc$  the $k$-rational point of $\mfc$ defined as the image of the zero element of $\mfg$ via $\chi : \mfg \migi \mfc$.
Note that $\mfc$ has a $\mbG_m$-action  arising from the homotheties on $\mfg$, and $[0]_\mfc$ is invariant under this action. 

We denote by  $(\Omega_{\mbP^{D_{{\bf x}}\text{-}\mr{log}}/k}^{\otimes p})^\times$   the $\mbG_m$-bundle on  $\mbP$ corresponding to the line bundle $\Omega_{\mbP^{D_{{\bf x}}\text{-}\mr{log}}/k}^{\otimes p}$.
Then, we obtain  the twist $(\Omega_{\mbP^{D_{{\bf x}}\text{-}\mr{log}}/k}^{\otimes p})^\times \times^{\mbG_m} \mfc$ of $\mfc$ by  $(\Omega_{\mbP^{D_{{\bf x}}\text{-}\mr{log}}/k}^{\otimes p})^\times$.
The point $[0]_\mfc$ of $\mfc$ determines a global section 
\begin{align}
[0]_{\mfc, \mbP^{D_{{\bf x}}\text{-}\mr{log}}} : \mbP \migi (\Omega_{\mbP^{D_{{\bf x}}\text{-}\mr{log}}/k}^{\otimes p})^\times \times^{\mbG_m} \mfc
\end{align}
 of the natural projection $(\Omega_{\mbP^{D_{{\bf x}}\text{-}\mr{log}}/k}^{\otimes p})^\times \times^{\mbG_m} \mfc \migi \mbP$.
For a Miura $G$-oper $\widehat{\mcE}^\spadesuit := (\mcE, \nabla, \mcE_B, \mcE'_B)$ on $\mbP^{D_{{\bf x}}\text{-}\mr{log}}$,
the $p$-curvature of $\nabla$ can be expressed as a global section $\psi$ of  $\Gamma (\mbP, \Omega_{\mbP^{D_{{\bf x}}\text{-}\mr{log}}/k}^{\otimes p} \otimes \mfg_\mcE)$;
it induces, via the quotient $\chi$, a global section 
\begin{align} \label{Ee100}
\psi^\chi : \mbP \migi (\Omega_{\mbP^{D_{{\bf x}}\text{-}\mr{log}}/k}^{\otimes p})^\times \times^{\mbG_m} \mfc.
\end{align}
We say that  $\widehat{\mcE}^\spadesuit$ is {\bf $p$-nilpotent}
(cf. ~\cite[Definition 3.8.1]{Wak7})  if 
$\psi^\chi = [0]_{\mfc, \mbP^{D_{{\bf x}}\text{-}\mr{log}}}$.
Thus,
we obtain the subset 
\begin{align} \label{Ee556}
G \text{-}\mr{MOp}(\mbP^{D_{{\bf x}}\text{-}\mr{log}}; - \check{\pmb{\lambda}}-\check{\rho})_{\mr{triv}, + \pmb{\alpha}^W}^{^{p\text{-}\mr{nilp}}}
\end{align}
of $G \text{-}\mr{MOp}(\mbP^{D_{{\bf x}}\text{-}\mr{log}}; - \check{\pmb{\lambda}}-\check{\rho})_{\mr{triv}, + \pmb{\alpha}^W}$ consisting of $p$-nilpotent Miura $G$-opers.

\SSP
\bpr \label{pgT011}
Suppose that the equality (\ref{W100}) holds.
Then, the   natural inclusion
\begin{align} \label{Lpdopa}
G \text{-}\mr{MOp}(\mbP^{D_{{\bf x}}\text{-}\mr{log}}; - \check{\pmb{\lambda}}-\check{\rho})_{\mr{triv}, + \pmb{\alpha}^W}^{^{p\text{-}\mr{nilp}}}
\migiincl G \text{-}\mr{MOp}(\mbP^{D_{{\bf x}}\text{-}\mr{log}}; - \check{\pmb{\lambda}}-\check{\rho})_{\mr{triv}, + \pmb{\alpha}^W}
\end{align}
is bijective.
In particular, 
the bijection (\ref{Ee55}) induces  a bijection 
 \begin{align}
G \text{-} \mr{BA}_{\check{\pmb{\lambda}}, \pmb{\alpha}} \isom 
G \text{-}\mr{MOp}(\mbP^{D_{{\bf x}}\text{-}\mr{log}}; - \check{\pmb{\lambda}}-\check{\rho})_{\mr{triv}, + \pmb{\alpha}^W}^{^{p\text{-}\mr{nilp}}}.
 \end{align}
  \epr
\begin{proof}
Let ${\bf z}:= (z_1, \cdots, z_m)$ be an element of 
$G \text{-} \mr{BA}_{\check{\pmb{\lambda}}, \pmb{\alpha}}$, and denote by
$\widehat{\mcE}^\spadesuit := (\mcE, \nabla, \mcE_B, \mcE'_B)$ the corresponding Miura $G$-oper via (\ref{Ee55}).
Denote by $B^{-}$  the opposite Borel subgroup of $B$ relative to $T$  and by  $\mfb^{-}$ its Lie algebra.
Then, the following square diagram is commutative:
\begin{align} \label{Ej698}
\begin{CD}
\mfb^{-}@> \mr{incl.} >> \mfg
\\
@VVV @VV \chi V
\\
\mft @>>> \mfc,
\end{CD}
\end{align}
where the left-hand vertical and lower horizontal arrows denote the natural quotients.
Since $B \cap B^{-} = T$, two reductions $\mcE_B$, $\mcE'_B$ induce a $T$-reduction $\mcE_T$ of $\mcE$.
The $B$-reduction $\mcE'_B$ is horizontal,  so 
the $p$-curvature $\psi$ of $\nabla$  lies in  $\Gamma (\mbP, \Omega_{\mbP^{D\text{-}\mr{log}}/k}^{\otimes p} \otimes \mfb^{-}_{\mcE_T})  \left(\subseteq \Gamma (\mbP, \Omega_{\mbP^{D\text{-}\mr{log}}/k}^{\otimes p} \otimes \mfg_\mcE) \right)$, where $D := D_{({\bf x}, {\bf z})}$.
Denote by $\overline{\psi}$ the image of $\psi$ in $\Gamma (\mbP, \Omega_{\mbP^{D\text{-}\mr{log}}/k}^{\otimes p} \otimes \mft_{\mcE_T})$ via the natural quotient $\mfb^{-} \migisurj \mft$.
Since the set-theoretic preimage of $[0]_\mfc$ via the quotient $\mft \migisurj \mfc$ consists   exactly of the zero element,
the commutativity of (\ref{Ej698}) implies that  the equality  $\overline{\psi} =0$ holds if and only if  $\widehat{\mcE}^\spadesuit$ is $p$-nilpotent.
Here, notice (cf. the proofs of Proposition  \ref{P05} and Theorem \ref{TT011}) that the element $\overline{\psi}$ (restricted to $\mbP \setminus \{ \infty \}$) may be identified with  the $p$-curvature of the connection $\nabla'$ on the trivial $T$-bundle expressed as
\begin{align}
\nabla' = \partial_x + f(x), \hspace{5mm} \text{where} \ f(x) :=  \sum_{i=1}^r \frac{- \check{\lambda}_i}{x -x_i} + \sum_{j=1}^m \frac{\check{\alpha}_j}{x- z_j}.
\end{align}
If  $\sigma$ denotes the base-change  $X^{(1)} \migi X$  via  the Frobenius automorphism of $k$ and $C_{\mbP^{D\text{-}\mr{log}}/k}$ denotes the Cartier operator of $\mbP^{D\text{-}\mr{log}}/k$ (cf. (\ref{eoaipfp})), then  the equality 
$C_{\mbP^{D\text{-}\mr{log}}/k} (f(x)) = \sigma^*(f(x))$ holds.
Hence,  by  ~\cite[Corollary 7.1.3]{K2},  $\nabla'$ has vanishing $p$-curvature.
It follows from the above discussion that $\widehat{\mcE}^\spadesuit$ is $p$-nilpotent.
This shows the bijectivity of 
the inclusion (\ref{Lpdopa}), so we finish 
 the proof of the proposition.
\end{proof}

\vspace{10mm}
\section{Tango structures} \SSP

In this section, we recall the notion of a Tango structure and consider a bijective correspondence between Tango structures and  solutions to the  Bethe ansatz equations for $\mfg= \mfs \mfl_2$.
Also, by means of this  correspondence, 
we construct  new examples of Tango structures. 
Hereinafter, we suppose  that $p\geq 3$.

\LSP
\subsection{Tango structures} \label{SS131}


Let $X$   be a connected 
proper 
 smooth curve over $k$ of genus $g_X \left(\geq 0\right)$.
 Denote by $X^{(1)}$ the Frobenius twist of $X$ relative to $k$ and by $F : X \migi X^{(1)}$ the relative Frobenius morphism of $X$.
 Also, denote by $\mcB_{X/k} \left(\subseteq \Omega_{X/k}\right)$ the sheaf of locally exact
   $1$-forms on $X$ relative to $k$.
The direct image  $F_* (\Omega_{X/k})$ of $\Omega_{X/k}$ forms a vector bundle on $X^{(1)}$ of rank $p$.
The sheaf  $\mcB_{X/k}$ may be considered, via the underlying homeomorphism of $F$, as a subbundle of $F_*(\Omega_{X/k})$ of rank $p-1$.

 Now, let $\mcL$ be a line subbundle of $\mcB_{X/k}$.
 Consider the $\mcO_{X^{(1)}}$-linear composite
 \begin{align}
 \mcL \migiincl F_* (\mcB_{X/k}) \migiincl F_*(\Omega_{X/k}),
 \end{align}
 where the first arrow denotes the natural inclusion and the second arrow denotes the morphism obtained by applying the functor $F_*(-)$ to the natural inclusion $\mcB_{X/k} \migiincl \Omega_{X/k}$.
This composite corresponds to a morphism
\begin{align}
\xi_\mcL : F^*(\mcL) \migi \Omega_{X/k}
\end{align}
via the adjunction relation  ``$F^*(-) \dashv F_*(-)$".

\SSP
\bde  \label{W400}
\begin{itemize}
\item[(i)]
We shall say that 
$\mcL$
  is a {\bf Tango structure} on $X$ if 
  $\xi_\mcL$ is an isomorphism.
\item[(ii)]
A {\bf Tango curve} is a connected proper smooth curve over $k$ admitting a Tango structure.
\end{itemize}
 \ede
\SSP

\begin{rema}\label{Rr0911}
 Let $\mcL$ be a line subbundle of $\mcB_{X/k}$.
 Since $\mr{deg}(\Omega_{X/k}) = 2g_X -2$, one may  verify that {\it $\mcL$ defines  a Tango structure on $X$ if and only if $\mcL$ has  degree $\frac{2g_X-2}{p}$}.
 In particular,   $X$ admits  no Tango structure when  $p \nmid g_X-1$.
   Y. Hoshi  proved that there exists a Tango curve of genus $g$ if and only if $p \mid g-1$  (cf. ~\cite[Theorem 1]{Hos}).
   \end{rema}
\SSP

Next, we shall recall  the definition of a pre-Tango structure in the sense of ~\cite[Definition 5.3.1]{Wak7}.
 Let $r$ be a nonnegative integer  and 
  ${\bf x} := (x_1, \cdots, x_r)$ 
   an ordered collection of distinct closed points of $X$.
(We take ${\bf x}:=\emptyset$ if $r=0$.)
If there is no fear of confusion, we shall write
$X^{\mr{log}}$ instead of  $X^{D_{{\bf x}}\text{-} \mr{log}}$, where $D_{{\bf x}} := \sum_{i=1}^r [x_i]$.
The data  ${\bf x}$ induces, via base-change,  a collection of closed points ${\bf x}^{(1)} := (x_1^{(1)}, \cdots, x_r^{(1)})$ in $X^{(1)}$, which determines a log structure on $X^{(1)}$; we denote  the resulting  log scheme by
$X^{(1) \mr{log}}$.

Moreover, let
\begin{align} \label{eoaipfp}
C_{X^\mr{log}/k} : F_{*} (\Omega_{X^{\mr{log}}/k}) \migi \Omega_{X^{(1)\mr{log}}/k}
\end{align}
be the {\it Cartier operator} of $X^\mr{log}/k$.
To be precise, $C_{X^\mr{log}/k}$ is a unique $\mcO_{X^{(1)}}$-linear  morphism 
whose composite with the injection
$\Omega_{X^{(1) \mr{log}}/k} \migiincl  \Omega_{X^{(1)\mr{log}}/k} \otimes F_*(\mcO_X)$
induced by the natural injection 
$\mcO_{X^{(1)}} \migiincl F_*(\mcO_X)$ coincides with the Cartier operator associated with the trivial flat bundle $(\mcO_X, d)$  in the sense of ~\cite[Proposition 1.2.4]{Og}.

\SSP
\bde \label{Ww551}
 A {\bf pre-Tango structure} on $(X, {\bf x})$
 is a $D_{{\bf x}}$-log connection $\nabla$ on $\Omega_{X^\mr{log}/k}$
   with vanishing $p$-curvature satisfying that
 $\mr{Ker}(\nabla) \subseteq \mr{Ker}(C_{X^\mr{log}/k})$.
  \ede
\SSP

\begin{rema}\label{Rr0922}
 According to ~\cite[Proposition 5.3.2]{Wak7}, 
{\it  pre-Tango structures on $(X, \emptyset)$ in the above sense  correspond bijectively to  Tango structures on $X$ in the sense of Definition \ref{W400}}.

 Indeed, let $\nabla$ be a pre-Tango structure on $(X, \emptyset)$.
 The $\mcO_{X^{(1)}}$-module $\mr{Ker}(\nabla)$ is 
 contained in $F_*(\mcB_{X/k})  \left(= F_*(\mr{Ker}(C_{X/k})) \right)$.
 Since $\nabla$ induces  an injection  
   $F_*(\Omega_{X/k})/\mr{Ker}(\nabla) \migiincl F_*(\Omega_{X/k}^{\otimes 2})$ (which implies that $F_*(\Omega_{X/k})/\mr{Ker}(\nabla)$ is a vector bundle),
  the  $\mcO_{X^{(1)}}$-submodule $F_*(\mcB_{X/k})/\mr{Ker}(\nabla)$ of $F_*(\Omega_{X/k})/\mr{Ker}(\nabla)$ turns out to be  a vector bundle.
  Namely,  $\mr{Ker}(\nabla)$ specifies a subbundle of $F_*(\mcB_{X/k})$.
  Moreover, 
    the condition that $\nabla$ has vanishing $p$-curvature implies that 
    the morphism 
  $\xi_{\mr{Ker}(\nabla)} : F^*(\mr{Ker}(\nabla)) \migi \Omega_{X/k}$ is an isomorphism.
 Hence, $\mr{Ker}(\nabla)$ defines  a Tango structure.
 
 Conversely, let  $\mcL$ be  a Tango structure on $X$.
 The line bundle $F^*(\mcL)$ has uniquely a  connection determined by the condition that the sections in $F^{-1}(\mcL)$ are horizontal.
 The connection $\nabla_\mcL$ corresponding, via the isomorphism $\xi_\mcL$, to this connection specifies a pre-Tango structure on $(X, \emptyset)$ (because of the equality $F_*(\mcB_{X/k}) = F_*(\mr{Ker}(C_{X/k}))$).
The resulting assignments $\nabla \mapsto \mr{Ker}(\nabla)$ and $\mcL \mapsto \nabla_\mcL$ together  give the desired  correspondence.
  \end{rema}
\SSP

Let us fix an ordered collection 
$\pmb{l} := (l_1, \cdots, l_r)\in \mbF_p^{r}$.
(We take  $\pmb{l}:=\emptyset$ if $r=0$.)
Denote by 
\begin{align} \label{Ww550}
\mr{Tang} (X; {\bf x}; \pmb{l})
\end{align}
the set of pre-Tango structures on $(X, {\bf x})$ 
  whose monodromy operator 
  at $x_i$  coincides with $l_i$ for every  $i \in \{1, \cdots, r \}$.

\SSP
\bpr \label{W31}
\begin{itemize}
\item[(i)]
 Let $m$ be a positive integer and let ${\bf z} := (z_1, \cdots, z_m) \in C^m ({\bf x})$.
 Then, there exists a canonical  bijection 
 \begin{align}
 \mr{Tang} (X; ({\bf x}, {\bf z}); (\pmb{l}, (-1)^{m})) \isom \mr{Tang} (X; {\bf x}; \pmb{l}).
 \end{align}
  \item[(ii)]
    There exists a bijective correspondence between  the set $\mr{Tang} (X; {\bf x}; (-1)^{r})$ and the set of Tango structures on $X$.
  \end{itemize}
 \epr
\begin{proof}
 First, we shall  consider assertion (i).
  Let  us take a pre-Tango structure  $\nabla$  classified by 
 $ \mr{Tang} (X; ({\bf x}, {\bf z}); (\pmb{l}, (-1)^{m}))$.
 Denote by $\nabla_{- {\bf z}}$ 
  the $D_{({\bf x}, {\bf z})}$-log connection on $\Omega_{X^{D_{{\bf x}}\text{-}\mr{log}}/k}$
    obtained by restricting $\nabla$ via the inclusion 
 $\Omega_{X^{D_{{\bf x}}\text{-}\mr{log}}/k} \migiincl \Omega_{X^{D_{({\bf x}, {\bf z})}\text{-}\mr{log}}/k}$.
The monodromy operator  of $\nabla_{- {\bf z}}$ at $z_j$ (for each $j \in \{1, \cdots, m \}$)
 coincides with $0$.
 Hence, 
 $\nabla_{- {\bf z}}$ may be thought of as a $D_{{\bf x}}$-log connection.
 Moreover, since the equality 
 \begin{align} \label{Www300}
 \mr{Ker}(C_{X^{D_{{\bf x}}\text{-}\mr{log}}/k}) = \mr{Ker}(C_{X^{D_{({\bf x}, {\bf z})}\text{-}\mr{log}}/k})
 \end{align}
 holds, 
  $\nabla_{- {\bf z}}$ specifies a pre-Tango structure classified by $\mr{Tang} (X; {\bf x}; \pmb{l})$.

 Conversely, let us take a pre-Tango structure $\nabla$ classified  by $\mr{Tang} (X; {\bf x}; \pmb{l})$;  
 it may be regarded as a $D_{({\bf x}, {\bf z})}$-log connection whose monodromy  operator at $z_j$  is equal to $0$ for every $j \in \{1, \cdots, m \}$.
 Then, there exists uniquely a $D_{({\bf x}, {\bf z})}$-log connection $\nabla_{+ {\bf z}}$ on 
 $\Omega_{X^{D_{({\bf x}, {\bf z})}\text{-}\mr{log}}/k}$ whose restriction to $\Omega_{X^{D_{{\bf x}}\text{-}\mr{log}}/k}$ coincides with $\nabla$.
 It follows  from the equality (\ref{Www300})  that $\nabla_{+ {\bf z}}$ specifies a pre-Tango structure classified by  $ \mr{Tang} (X; ({\bf x}, {\bf z}); (\pmb{l}, (-1)^{m}))$.
 
 The resulting assignments $\nabla \mapsto  \nabla_{- {\bf z}}$ and  $\nabla \mapsto \nabla_{+ {\bf z}}$   give the desired bijection.
 This completes the proof of assertion (i).
 
 Assertion (ii) follows from assertion (i) of the case where $m =0$ together with the result mentioned in Remark \ref{Rr0922}.
\end{proof}
\SSP

The following assertion is  an important  property of Tango structures, which provides a relationship with dormant generic Miura $\mr{PGL}_2$-opers.

\SSP
\bt[cf.  ~\cite{Wak7}, Theorem A]\label{Fg011} 
Let us keep the above notation, and suppose that  $\pmb{l} 
\in (\mbF_p^\times)^{r}$.
 Then, 
there exists a canonical bijection
\begin{align} \label{W49}
\mr{Tang} (X; {\bf x}; \pmb{l}) \isom \mr{PGL}_2 \text{-} \mr{MOp} (X^{D_{{\bf x}}\text{-} \mr{log}}; \pmb{l} \check{\rho})^{^\mr{Zzz...}}_{\mr{gen}},
\end{align}
where $\check{\rho} = \begin{pmatrix} \frac{1}{2} & 0 \\ 0 & -\frac{1}{2}\end{pmatrix}$ and  $\pmb{l} \check{\rho} := (l_1   \check{\rho}, \cdots, l_r  \check{\rho}) \in (\mft_\mr{reg}^F)^{r}$.
  \et

\LSP
\subsection{Pull-back via tamely ramified coverings} \label{SS400}

We  consider  the pull-back of a pre-Tango structure by  a tamely ramified covering of $X$.
Suppose that we are given a collection of data
\begin{align} \label{ww1}
(Y, {\bf y}, \pi),
\end{align}
where 
\begin{itemize}
\item
  $Y$ denotes another connected proper smooth curve over $k$;
  \item
  ${\bf y}:= (y_1, \cdots, y_s)$ ($s >1$) denotes  an ordered collection of distinct closed points of $Y$;
  \item
  $ \pi : Y \migi X$ denotes a {\it tamely ramified} covering such that $\pi^{-1}(\bigcup_{i=1}^r \{ x_i \}) = \bigcup_{j=1}^s \{ y_j \}$ and $\pi$ is \'{e}tale away from $\bigcup_{i=1}^r \{ x_i \}$.
 \end{itemize}
Denote by $q : \{1, \cdots, s \} \migi \{ 1,\cdots, r \}$ the map  determined by 
$\pi (y_j) = x_{q(j)}$ ($j =1, \cdots, s$).
For each $j \in \{1, \cdots, s \}$, denote by $R_{j}$ the ramification index of $\pi$ at $y_j$ (hence $p \nmid R_j$).
The morphism 
 $\pi$ extends to a log \'{e}tale morphism
 $\pi^\mr{log} : Y^{D_{{\bf y}}\text{-} \mr{log}} \migi X^{D_{{\bf x}}\text{-} \mr{log}}$, and hence,  
 the natural morphism $\pi^*(\Omega_{X^{D_{{\bf x}}\text{-} \mr{log}}/k}) \migi \Omega_{Y^{D_{{\bf y}}\text{-} \mr{log}}/k}$ is an isomorphism.

 Next, let $\nabla$ be an element of $\mr{Tang} (X; {\bf x}; \pmb{l})$.
 The pull-back
    of $\nabla$ via $\pi$ determines, under the isomorphism $\pi^*(\Omega_{X^{D_{{\bf x}}\text{-} \mr{log}}/k}) \isom  \Omega_{Y^{D_{{\bf y}}\text{-} \mr{log}}/k}$, 
  a $D_{{\bf y}}$-log connection 
  \begin{align}
  \pi^*(\nabla) : \Omega_{Y^{D_{{\bf y}}\text{-} \mr{log}}/k} \migi \Omega_{Y^{D_{{\bf y}}\text{-} \mr{log}}/k}^{\otimes 2}
  \end{align}
   on $\Omega_{Y^{D_{{\bf y}}\text{-} \mr{log}}/k}$.
Here, observe that  both  the relative Frobenius morphisms and  the Cartier operators  of $X \setminus \bigcup_{i=1}^r \{ x_i \}$ and $Y \setminus \bigcup_{i=1}^s \{ y_i \}$ are compatible (in an evident sense) via $\pi$.
This implies that
  $\mr{Ker}(\pi^*(\nabla)) \subseteq \mr{Ker}(C_{Y^{D_{{\bf y}}\text{-} \mr{log}}/k})$, and hence, 
 $\pi^*(\nabla)$ specifies a pre-Tango structure on $(Y, {\bf y})$.

\SSP
\bpr \label{ww2}
Write $\pmb{l}' := (l_{q(1)}R_{1}, \cdots, l_{q(s)} R_s) \in \mbF_p^{s}$.
 Then, the  assignment $\nabla \mapsto \pi^*(\nabla)$ discussed above defines an injection
 \begin{align}
 \mr{Tang} (X; {\bf x}; \pmb{l}) \migiincl  \mr{Tang} (Y; {\bf y}; \pmb{l}'). 
 \end{align}
  \epr
\begin{proof}
Since the injectivity may be immediately verified, it suffices to prove that the monodromy operator of the pre-Tango structure $\pi^*(\nabla)$ at $y_j$ (where $\nabla \in  \mr{Tang} (X; {\bf x}; \pmb{l})$, $j \in \{1, \cdots, s \}$) is $l_{q(j)}R_j$.
For simplicity, we write $x := x_{q(j)}$, $y := y_j$, $R := R_j$, and $l := l_{q(j)}$.
The formally local description of  $\pi$ at $y$ may be given by $\mbD_{y_{}} := \mr{Spec}(k[\![t^{1/R}]\!]) \migi \mr{Spec}(k[\![t]\!]) := \mbD_{x_{}}$ corresponding to the natural inclusion $k[\![t]\!] \migiincl k[\![t^{1/R}]\!]$, where $t$ denotes  a  formal coordinate in $X$  at $x_{}$.
By a suitable trivialization $\Omega_{X^{D_{{\bf x}}\text{-} \mr{log}}/k}  |_{\mbD_{x_{}}} \cong \mcO_{\mbD_{x_{}}}$,
the restriction $\nabla |_{\mbD_{x_{}}}$ of $\nabla$ to $\mbD_{x_{}}$ 
 is  expressed  as $\nabla = d +  l \cdot  \frac{dt}{t}$ (cf. ~\cite[Corollary 2.10]{Os}).
  Since
  \begin{align}
  \frac{d t}{t} = \frac{d (t^{1/R_{}})^R}{ (t^{1/R})^R} = \frac{R \cdot (t^{1/R})^{R-1} \cdot  d t^{1/R}}{(t^{1/R})^R} = R\cdot \frac{d t^{1/R}}{t^{1/R}}, 
  \end{align}
 the restriction $\pi^*(\nabla)|_{\mbD_{y_{}}}$ of $\pi^*(\nabla)$ to  $\mbD_{y_{}}$ satisfies  $\pi^*(\nabla)|_{\mbD_{y_{}}} = d + l R \cdot \frac{dt^{1/R}}{t^{1/R}}$.
It follows  that the monodromy operator of $\pi^*(\nabla)$  at $y$ coincides with $lR$.
This  completes the proof of the assertion.
\end{proof}

\SSP
\begin{exa}\label{Er09}
 We focus on  the case where $X =\mbP$.
 Let $\pmb{l} := (l_1, \cdots, l_{r+1})$ ($r\geq 0$) be an element of $\mbF^{r+1}_p$ and     ${\bf x} := (x_1, \cdots, x_{r+1})$  an ordered collection of distinct $r+1$ closed points of $\mbP$ with $x_{r+1} := \infty$.
 Let us consider the desingularization $Y$ of the plane curve defined by
 \begin{align}
 y^n = a \cdot \prod_{i=1}^r (x - x_i),
  \end{align}
 where $a \in k^{\times}$, $p \nmid n$,   and $(x, y)$ are an inhomogeneous coordinate of the projective plane $\mbP^2$.
Denote by $\pi : Y \migi \mbP$ the projection given by $(x, y) \mapsto x$.
For each  $i\in \{1, \cdots, r \}$, we  define $y_i$ to be  the unique point of $Y$ lying over $x_i$, i.e., $\pi (y_i) =x_i$.
(The ramification index of $\pi$ at $y_i$  coincides with  $n$.)
Also, let $y_{r+1}, \cdots, y_{r+r'}$ ($r'  \geq 1$) be the set of points of $Y$ lying  over $x_{r+1} \ (=\infty)$. 
Write ${\bf y} := (y_1, \cdots, y_{r+r'})$.
According to ~\cite[Proposition 6.3.1]{Sti},  the ramification index of $\pi$ at $y_{r + i'}$ (for each $i' \in \{1, \cdots, r'\}$) is $\frac{n}{d}$, where $d := \mr{gcd} (n, r)$.
Moreover, the genus of $Y$ is given by $\frac{(n-1)(r-1)}{2} - \frac{\mr{gcd}(n, r)-1}{2}$.

Now, let us fix $\pmb{l} := (l_1, \cdots, l_{r+1}) \in \mbF_p^{r+1}$ and write
\begin{align}
\pmb{l}' := (l_1 n, l_2 n, \cdots, l_r n, \frac{l_{r+1}n}{d}, \cdots, \frac{l_{r+1}n}{d}) \in \mbF_p^{r +r'}
\end{align}
 (where each of  the last $r'$ factors is $\frac{l_{r+1}n}{d}$).
Then, it follows from Proposition \ref{ww2} that  
 the assignment $\nabla \mapsto \pi^*(\nabla)$ define an injection
\begin{align}
\mr{Tang} (\mbP; {\bf x}; \pmb{l}) \migiincl \mr{Tang} (Y; {\bf y}; \pmb{l}').
\end{align}
In particular, $\mr{Tang} (Y; {\bf y}; \pmb{l}')$ is nonempty unless   $\mr{Tang} (\mbP; {\bf x}; \pmb{l})$ is empty.
 \end{exa}

\LSP
\subsection{Examples of Tango structures} \label{Ss131}

Let us  show that
 the  solutions of the Bethe ansatz equations of a certain type come from
    dormant Miura $\mr{PGL}_2$-opers.

To this end, we apply the bijective correspondence (\ref{Ee990}) in the case of $G=\mr{PGL}_2$ and  $p\geq 3$ 
 (i.e.,   the condition   $(*)_{\mr{PGL}_2}$ is satisfied).
Recall that  
 $\check{\rho} = \begin{pmatrix} \frac{1}{2} & 0 \\ 0 & -\frac{1}{2}\end{pmatrix}$,  $\Gamma = \{ \alpha \}$, and  $\check{\alpha} = 2 \check{\rho} = \begin{pmatrix} 1 & 0 \\ 0 & -1\end{pmatrix}$.
If $r=0$, then the equality (\ref{W100}) holds if and only if $m\equiv 1$ mod $p$.
Now, let us take a positive integer $l$ (and consider the case of $m= lp+1$).
We  combine   Theorems \ref{TT011},  \ref{Fg011}, and Proposition \ref{W31}, (i), to obtain  the following composite bijection
\begin{align} \label{Ee860}
\mr{PGL}_2\text{-} \mr{BA}_{0, \alpha^{lp+1}} &\isom
\coprod_{{\bf z} \in C^{lp+1}(\infty)}\mr{PGL}_2 \text{-}\mr{MOp}(\mbP^{D_{(\infty, {\bf z})}}; (-\check{\rho}, \check{\rho}^{lp+1}))
\\
&\isom 
\coprod_{{\bf z} \in C^{lp+1}(\infty)}\mr{Tang} (\mbP; (\infty, {\bf z});  (-1, 1^{lp+1})) \notag  \\
&\isom 
\coprod_{{\bf z} \in C^{lp+1}(\infty)}\mr{Tang} (\mbP;  {\bf z};  1^{lp+1}). \notag 
\end{align}
In particular,  it follows from the discussion in Remark  \ref{Rr319} that
 $(\mbP, {\bf z})$ admits  a pre-Tango structure  whose  monodromy operators are given by  $1^{lp+1}$ if and only if 
  the equality $f''(x) =0$ holds, where $f (x) := \prod_{i=1}^{lp +1} (x - z_i)$.

Next, let $a$, $b$ be  positive integers with $\mr{gcd}(a, b p -1)=1$
and $h (x)  \left(\in k[x] \right)$
  a monic polynomial  of degree $ap$ with
 $\mr{gcd}(h (x), h'(x))=1$.
Denote by
\begin{align}
Y
\end{align}
 the smooth   curve defined by
    the equation 
 \begin{align}
 y^{b p-1} = h (x).
 \end{align}
(If $a > b$, then the point at infinity   is singular, and hence, we  need to replace 
this  plane  curve 
 by its desingularization to obtain a smooth curve.)
Since $\mr{gcd}(ap, b p -1) = \mr{gcd}(a, bp-1)=1$,  there is  only one point $\infty_Y$ at infinity in $Y$ (cf. ~\cite[Proposition 6.3.1]{Sti}).
Here, let us take 
 an ordered collection $(z_1, \cdots, z_{ap})$  of elements of $k$ (i.e., closed points in $\mr{Spec}(k[x]) = \mbP \setminus \{ \infty \}$) 
 satisfying  the equality $h (x ) =\prod_{i=1}^{a p}(x - z_i)$.
The assumption $\mr{gcd}(h(x), h'(x))=1$ implies that $z_{i} \neq z_{i'}$ if $i \neq i'$.
Denote by $\pi :  Y\migi \mbP$ the natural projection $(x, y) \mapsto x$, which is tamely ramified.
Also,  denote   by ${\bf y} := (y_1, \cdots, y_{a p+1})$, where ${\bf z}:= (z_1, \cdots, z_{ap+1})$ and $z_{ap+1} := \infty$, the ordered collection of distinct points in $Y$ determined uniquely by $\pi (y_i) = z_i$ for any $i =1, \cdots, ap +1$.
Moreover, denote by $\nabla$ the $D_{{\bf z}}$-log connection on $\Omega_{\mbP^{D_{{\bf z}} \text{-} \mr{log}}/k}$ defined as  
\begin{align}
\nabla = \partial_x + \sum_{i=1}^{a p} \frac{1}{x - z_i}
\end{align}
 under the identification $\Omega_{\mbP^{D_{{\bf z}} \text{-} \mr{log}}/k} |_{\mbP \setminus \{ \infty \}}\isom \mcO_{\mbP \setminus \{ \infty \}}$ given by $d x \mapsto 1$.
By passing to  the isomorphism $\pi^*(\Omega_{\mbP^{D_{{\bf z}} \text{-} \mr{log}}/k}) \isom \Omega_{Y^{D_{{\bf y}} \text{-} \mr{log}}/k}$ induced by $\pi$,
we obtain a $D_{{\bf y}}$-log connection 
\begin{align}
\pi^*(\nabla) : \Omega_{Y^{D_{{\bf y}} \text{-} \mr{log}}/k} \migi \Omega_{Y^{D_{{\bf y}} \text{-} \mr{log}}/k}^{\otimes 2}
\end{align}
 on $\Omega_{Y^{D_{{\bf y}} \text{-} \mr{log}}/k}$
defined to be the pull-back of $\nabla$.
In particular, we obtain an $\mcO_{Y^{(1)}}$-submodule 
\begin{align}
\mr{Ker}(\pi^*(\nabla))
\end{align}
 of $F_{Y/k*}(\Omega_{Y^{D_{{\bf y}} \text{-} \mr{log}}/k})$, where $F_{Y/k} : Y \migi Y^{(1)}$ denotes the relative Frobenius morphism  of $Y$ relative to $k$.
Then, the following assertion, i.e., Theorem \ref{ThB},  holds.

\SSP
\bt  \label{T041} 
 Suppose that $h''(x) =0$, or equivalently, $(z_1, \cdots, z_{ap})$ specifies an element of $\mr{PGL}_2\text{-}\mr{BA}_{0, \alpha^{ap}}$.
 Then, $\mr{Ker}(\pi^*(\nabla))$ forms a Tango structure on $Y$.
 In particular, $Y$ is a Tango curve.
 \et
\begin{proof}
Let us fix an element $\gamma \in k$ with $h (\gamma) \neq 0$. 
Consider the automorphism $\iota$ of $\mbP$ given by $x \mapsto \frac{1}{x-\gamma}$.
Then, $\iota ({\bf z}):=(\iota (z_1), \cdots, \iota (z_{ap+1})))$ are a collection of  distinct points in $\mbP \setminus \{ \infty \}$.
Since $f(x) := \prod_{i=1}^{ap+1} (x - \iota (z_i))$ coincides with $h (\gamma)^{-1} \cdot h (\frac{1}{x}+\gamma) \cdot  x^{ap+1}$,
we see that $\mr{gcd}(f(x), f'(x))=1$ and $f''(x) =0$.
It follows from Theorem  \ref{ThgA}
that
the $D_{\iota ({\bf z})}$-log connection on $\Omega_{\mbP^{D_{\iota({\bf z})} \text{-} \mr{log}}/k}$ determined  by $\partial_x + \sum_{i=1}^{ap+1} \frac{1}{x - \iota (z_i)}$ forms a pre-Tango structure on $(\mbP, \iota ({\bf z}))$.
By pulling-back via  $\iota$,
we obtain a pre-Tango structure on $(\mbP, {\bf z})$
 whose monodromy operators are   $1^{ap +1}$;
by construction,  it coincides with $\nabla$.
Hence, (since $bp-1 \equiv -1$ mod  $p$)
the pull-back $\pi^*(\nabla)$ forms a pre-Tango structure on $(Z, {\bf z})$ whose   monodromy operators are given by  $(-1)^{ap+1}$ (cf. Proposition \ref{ww2}).
It follows from  Proposition 
\ref{W31} (and its proof) that 
$\mr{Ker}(\pi^*(\nabla))$  specifies a Tango structure.
\end{proof}

\SSP
\begin{rema}\label{Rr09}
A well-known example of a  Tango structures can be found in some literature (cf., e.g., ~\cite[Example]{Ray} and  ~\cite[Example 1.3]{Muk}).
 This is constructed as follows.
 Let $l$ be an integer with $lp \geq 4$.
Also, let $f (x)$ be a polynomial of degree $l$ in one variable $x$ and let $Y$ 
 be the  plane curve defined by 
 \begin{align}
 y^{lp-1} = f(x^p) -x.
 \end{align}
One verifies  that $Y$ is a smooth curve having only one point $\infty$ at infinity and  $\Omega_{Y/k} = \mcO_Y ( lp (lp -3) \cdot \infty)$.
Then, the base-change  of the line bundle $\mcO (l(lp-3) \cdot \infty)$ via the absolute Frobenius morphism of $\mr{Spec}(k)$
turns out to specify a Tango  structure. 

Notice that $(f(x^p) -x)'=-1$ (which implies  $\mr{gcd}(f(x^p) -x, (f(x^p) -x)')=1$) and $(f(x^p) -x)''=0$.
Hence, 
 $Y$ is    
a specific type of 
 Tango curves  constructed in the discussion preceding Theorem \ref{T041}.
On the other hand, the curve defined, e.g.,  by the equation $y^{2p-1} = x^{2p}+ x^{p+1}+ ax^p +b x +c$ with $c \neq b (a-b)$ is
a Tango curve (by Theorem \ref{T041}), which gives 
 a new example.
Thus, our result allows us to obtain   infinitely many new (explicit!) examples of Tango curves.
  \end{rema}

\LSP
\subsection*{Acknowledgements} 
We would  like to thank
Professor Michel Gros for his useful comments.
We are very grateful for the many helpful conversations we had with 
 dormant Miura opers and Tango curves, living in the world of mathematics.
Our work was   partially  supported by the Grant-in-Aid for Scientific Research (KAKENHI No.\,21K13770).

\end{document}